\documentclass[12pt]{amsart}
\usepackage[indentafter]{titlesec}
\titleformat{name=\section}{}{\thetitle.}{0.8em}{\centering\scshape}
\titleformat{name=\subsection}[runin]{}{\thetitle.}{0.5em}{\bfseries}[.]
\titleformat{name=\subsubsection}[runin]{}{\thetitle.}{0.5em}{\itshape}[.]
\titleformat{name=\paragraph,numberless}[runin]{}{}{0em}{}[.]
\titlespacing{\paragraph}{0em}{0em}{0.5em}
\titleformat{name=\subparagraph,numberless}[runin]{}{}{0em}{}[.]
\titlespacing{\subparagraph}{0em}{0em}{0.5em}

\usepackage[
backend=biber,
style=alphabetic,
sorting=ynt
]{biblatex}
\usepackage[utf8]{inputenc}
\usepackage{amsmath,amssymb,amsthm}
\usepackage[margin=1in]{geometry}
\usepackage{mathbbol}
\usepackage{graphicx}
\usepackage{faktor}
\usepackage{stmaryrd}
\usepackage{tikz-cd}
\usepackage{makecell}
\usepackage{comment}

\usepackage[margin=1in]{geometry}

\usepackage[pdftex, plainpages=false, pdfpagelabels]{hyperref}
\hypersetup{
    linktocpage=true,
    colorlinks=true,
    bookmarks=true,
    citecolor=blue,
    urlcolor=blue,
    linkcolor=blue,
    citebordercolor={1 0 0},
    urlbordercolor={1 0 0},
    linkbordercolor={.7 .8 .8},
    breaklinks=true,
    pdfpagelabels=true,
    }

\usepackage{graphicx} 
\usepackage{amsmath} 
\usepackage{amsxtra} 
\usepackage{amssymb} 
\usepackage{amsthm} 
\usepackage{latexsym} 
\usepackage{setspace} 
\usepackage{nomencl} 
\usepackage[margin=1in]{geometry} 
\usepackage[titles]{tocloft} 
\newcommand{\C}{\mathbb{C}}
\newcommand{\PP}{\mathbb{P}}
\newcommand{\N}{\mathbb{N}}
\newcommand{\A}{\mathbb{A}}
\newcommand{\R}{\mathbb{R}}

\newcommand{\Z}{\mathbb{Z}}

\begin{document}

\newpage


\thispagestyle{empty}


\begin{center}

   \vspace{1cm}

\title{A Resolution of the Diagonal for Toric D-M Stacks}
\author{Reginald Anderson}
\email{rcanders@ksu.edu}
\maketitle


 \vspace{0.3cm}

\end{center}

\section{Abstract}

Beilinson gave a resolution of the diagonal for complex projective spaces which Bayer-Popescu-Sturmfels generalized to what they refer to as unimodular toric varieties. Here, we generalize Bayer-Popescu-Sturmfels resolution of the diagonal to the case of smooth toric varieties and global quotients of smooth toric varieties by a finite abelian group, considered as a toric Deligne-Mumford stack.


\section{Introduction}

In this paper, we study bounded derived categories of coherent sheaves $D^b_{Coh}(X)$ on a smooth projective variety $X$. The derived category $D^b_{Coh}(X)$ forms part of the B-side of Homological Mirror Symmetry (HMS). Beilinson's resolution of the diagonal \cite{Beilinson1978} was crucial in understanding $D^b(\PP^n)$. Since $\PP^n$ is a toric variety, we can ask for a generalization of Beilinson's resolution of the diagonal to any toric variety. However, this will be too broad, so we will restrict our attention to simplicial toric varieties. Bayer-Popescu-Sturmfels generalized Beilinson's resolution of the diagonal for unimodular toric varieties. We generalize Bayer-Popescu-Sturmfels' resolution of the diagonal to toric D-M stacks $\mathcal{X}(\boldsymbol{\Sigma}') = \left[ \faktor{X_\Sigma}{\mu} \right]$ associated to the global quotient of a smooth toric variety $X_\Sigma$ by a finite abelian group $\mu$. By deforming the resolution $\mathcal{F}^\bullet_{\mathcal{H}_L/L}$ of a unimouldar toric variety in the case $X_\Sigma$ is smooth and non-unimodular, we  can resolve the diagonal in families corresponding to a parameter $\epsilon$ in the effective cone of $X_\Sigma$. Since smooth toric varieties admit an open cover by affine spaces $\mathbb{A}^m$ as a variety, local diagonal objects for $\left[ \faktor{X_\Sigma}{\mu}\right]$ glue together to give a global diagonal object for $\mathcal{X}(\boldsymbol{\Sigma}')$.

\subsection{Results} 

This paper focuses on the following two theorems.

\textbf{Theorem}

\begin{enumerate}
\item[(A)] The complex $(\mathcal{F}^\bullet_{\mathcal{H}_L^\epsilon/L}, \partial^\epsilon)$ of $S$-free modules gives a resolution of the diagonal object $\pi(M_{\Lambda(L)})$ of $D^b_{Coh}(X_\Sigma \times X_\Sigma)$ for $X_\Sigma$ a smooth toric variety.\\
\item[(B)] The complex $(\mathcal{F}^\bullet_{\mathcal{H}_L^\epsilon/\tilde{L}}, \partial^\epsilon)$ resolves the diagonal object $\pi(M_{\Lambda(L)}) = M_{\Lambda(L)} \otimes_{S[\Lambda(\tilde{L})]} S$ of $D^b_{Coh}(X_\Sigma' \times X_\Sigma')$ for a simplicial toric variety $X_\Sigma' = \faktor{X_\Sigma}{\mu}$ viewed as a global quotient toric Deligne-Mumford stack where $X_\Sigma$ is a smooth toric variety and $\mu$ is a finite abelian group.
\end{enumerate}

\noindent 

Serre's theorem for projective space gives an equivalence between quasi-finitely generated graded modules over the homogeneous coordinate ring modulo the equivalence that $M \sim M'$ if there is an integer $d$ such that $M_{\geq d} \cong M'_{\geq d}$ \cite{hartshorne}. Cox's theorem gives an equivalence between graded modules over the homogeneous coordinate ring of a simplicial toric variety $X_\Sigma$, modulo torsion, and coherent sheaves on $X_\Sigma$ \cite{coxhgscoordinatering}. Specifically, Cox's theorem on the relation between coherent sheaves on a simplicial toric variety $X_\Sigma$ and graded modules over the homogeneous coordinate ring, modulo torsion states:

\begin{enumerate}
\item On a simplicial toric variety $X_\Sigma$, every coherent sheaf is of the form $\tilde{M}$ for some finitely generated module over the homogeneous coordinate ring $S$.
\item If $X_\Sigma$ is a smooth toric variety and $F$ is a finitely generated graded $S$-module, then $\tilde{F}=0$ iff there is some $k>0$ with $I_{irr}^k F = 0$ for $I_{irr}$ the irrelevant ideal of $X_\Sigma$.  
\end{enumerate}

We rely on part 2 of Cox's theorem to resolve the diagonal for $X_\Sigma$ smooth. \\

\noindent Cox's theorem has the following generalization for simplicial toric D-M stacks. Let the simplicial toric variety $X_\Sigma = \faktor{\left(\A^n \setminus V(I_{irr}) \right)}{ T}$ as a GIT quotient, with $\chi$ the associated toric D-M stack. Let $\mathcal{U}$ be the quasi-affine variety $\A^n \setminus V(I_{irr})$.\\

\begin{enumerate}
    \item There is an equivalence of categories between finitely generated modules over the homogeneous coordinate ring of $\mathcal{U}$ (modulo torsion) and coherent sheaves on $\mathcal{U}$.\\
    \item The category of graded, finitely generated modules over the homogeneous coordinate ring of $\mathcal{U}$ (modulo torsion) is equivalent to the category of equivariant coherent sheaves on $\mathcal{U}$.\\
    \item The category of graded, finitely generated modules over the homogeneous coordinate ring (mod torsion) of $\mathcal{U}$ is equivalent to coherent sheaves on the toric D-M stack $\chi$. 
    
\end{enumerate}

If $X_\Sigma$ is a simplicial toric variety, then we view it as a simplicial toric D-M stack. \\[1cm]

\noindent  Bayer-Popescu-Sturmfels\cite{bayer-popescu-sturmfels} give a cellular resolution of the Lawrence ideal $J_L$ corresponding to a unimodular lattice $L$. Imposing the requirement that $L$, the lattice of principal divisors of the projective toric variety be unimodular is fairly restrictive: In this case $L$ is given by the image of $B$ in the fundamental exact sequence\cite{C-L-S} \[ 0 \rightarrow M \stackrel{B}{\rightarrow} \Z^{|\Sigma(1)|} \stackrel{\pi}{\rightarrow} Cl(X_\Sigma) \rightarrow 0 \] 

\noindent 
satisfies the property that $B$ has linearly independent columns and every maximal minor of $B$ lies in the set $\{0,+1,-1\}$ (this is the sense of ``unimodularity" used in Bayer-Popescu-Sturmfels\cite{bayer-popescu-sturmfels}). This forms a proper subclass of smooth projective toric varieties, as $\PP^2$ blown up at a point and then blown up along the exceptional divisor is no longer unimodular. Also, any non-trivially weighted projective space is no longer unimodular. Throughout, the word ``unimodular" refers to this restricted sense of Bayer-Popescu-Sturmfels. \\

\noindent For a unimodular toric variety $X_\Sigma$,  the known resolution of $J_L$ gives a finite and minimal resolution of the ideal sheaf $\mathcal{O}_\Delta$ of the diagonal $\Delta \subset X_\Sigma \times X_\Sigma$ by sums of line bundles 
on $X_\Sigma\times X_\Sigma$. In the well-known example for $X_\Sigma=\PP^n_\C$, Beilinson's resolution gives the kernel of a Fourier-Mukai transform inducing the identity on $D^b_{Coh}(X_\Sigma)$. In this paper, we generalize this resolution to simplicial toric varieties $X_\Sigma'$ which are a global quotient of a smooth toric variety $X_\Sigma$ by a finite abelian group $\mu$ by viewing $X_\Sigma'$ as a toric D-M stack in Theorems (A) and (B) above.


For Theorem (B), the cofinite sublattice $\tilde{L}\subset L$ is chosen so that $\faktor{L}{\tilde{L}}\cong \mu$. Now the technique and deformations of the resolution of $J_L$ in Bayer-Popescu-Sturmfels generalize to give a resolution of the diagonal \[ (\mathcal{F}_{\mathcal{H}_L^\epsilon/\tilde{L}}, \partial) \] 

The proof is given in Section~\ref{section3.8} and culminates in Section~\ref{theorem3.8.3}.

\section{Known resolutions: minimal resolution of unimodular Lawrence ideal $J_L$ (Bayer-Popescu-Sturmfels)}
\label{sec:3.3}

Here, we recall the conventions of Bayer-Sturmfels\cite{bayer-sturmfels} and Bayer-Popescu-Sturmfels\cite{bayer-popescu-sturmfels} to resolve the diagonal for a unimodular toric variety. For $\Sigma$ a complete fan in the lattice $\Z^m$ and $X_\Sigma$ the associated complete normal unimodular toric variety with $\{\textbf{b}_1, \dots, \textbf{b}_n\} \subset \Sigma(1)$ the primitive ray generators of $\Sigma$, and $B$ the $n\times m$ matrix with row vectors $\textbf{b}_i$, each $\textbf{b}_i$ determines a torus-invariant Weil divisor $D_i$ on $X_i$, and the group $Cl(X)$ of torus-invariant Weil divisors modulo linear equivalence is given by\cite{bayer-popescu-sturmfels}

\begin{align} 0 \rightarrow \Z^m \stackrel{B}{\rightarrow} \Z^n \stackrel{\pi}{\rightarrow} Cl(X) \rightarrow 0, \end{align} \label{seq:toric.fund.seq.unimod}

where $\pi$ takes the $i$-th standard basis vector of $\Z^n$ to the linear equivalence class $[D_i]$ of the corresponding divisor $D_i$. \\

\textbf{Definiton}: The \underline{Lawrence ideal}\cite{bayer-popescu-sturmfels} $J_L$ corresponding to $L$, a sublattice of $\Z^n$ is 

\[ J_L = \left<\textbf{x}^\textbf{a}\textbf{y}^\textbf{b} - \textbf{x}^\textbf{b}\textbf{y}^\textbf{a} \text{ }|\text{ }\textbf{a}-\textbf{b} \in L \right> \subset S = k[x_1, \dots, x_n, y_1, \dots, y_n ] \]

for $k$ a field. Here $\textbf{x}^\textbf{a} = x_1^{a_1}x_2^{a_2}\cdots x_n^{a_n}$ for $\textbf{a}=(a_1, \dots, a_n) \in \N^n$.   

For $X_\Sigma$ unimodular with homogeneous coordinate ring $R = \C[x_\rho\text{ }|\text{ }\rho\in \Sigma(1)] \cong \C[x_1, \dots, x_n]$, the toric variety $X_\Sigma \times X_\Sigma$ has homogeneous coordinate ring 

\[ S = R \otimes_k R = k[x_1, \dots, x_n, y_1, \dots, y_n]. \] 

The diagonal embedding $X_\Sigma \subset X_\Sigma \times X_\Sigma$ defines a closed subscheme, and is represented by a $Cl(X) \times Cl(X)$-graded ideal $I_X$ in $S$. Here, $I_X = \ker(\psi)$ for $\psi: S = R\otimes_k R \rightarrow k[Cl(X)]\otimes R$ given by $\textbf{x}^{\textbf{u}}\textbf{y}^{\textbf{v}} = \textbf{x}^{\textbf{u}} \otimes \textbf{x}^{\textbf{v}} \mapsto [\textbf{u}]\otimes \textbf{x}^{\textbf{u}+\textbf{v}}$. That is,

\[ I_X = \left< \textbf{x}^\textbf{u}\textbf{y}^{\textbf{v}} - \textbf{x}^\textbf{v} \textbf{y}^\textbf{u} \text{ }|\text{ } \pi(\textbf{u}) = \pi(\textbf{v}) \text{ in } Cl(X) \right> \subset S. \]

\begin{proof}
First, we have the following \\

\textbf{Lemma}: Suppose $X = Spec(R)$ is an affine variety. Then the diagonal mapping $\Delta: V \rightarrow V \times V$ corresponds to the $\C$-algebra homomorphism $R \otimes R \rightarrow R$ given by $\sum r_1 \otimes r_2 \mapsto \sum r_1 r_2$, from the universal property of $X \times X$. \\

Next, locally, for $\Sigma \subset N_\R$, we have that the diagonal map $U_i \stackrel{\Delta}{\hookrightarrow} U_i \times U_i$ corresponds to 

\begin{align*}
\C[\sigma_i^\vee \cap M] \otimes_k \C[\sigma_i^\vee \cap M] \rightarrow \C [\sigma_i^\vee \cap M ] \text{ given by} \\
\sum r_1 \otimes r_2 \mapsto \sum r_1 r_2 . \end{align*}

A choice of $\sigma\in \Sigma$ gives\cite{C-L-S} the monomial \[ x^{\hat{\sigma} } = \prod_{\rho\not\in\sigma(1) } x_\rho \in S \] for $S$ the homogeneous coordinate ring\cite{coxhgscoordinatering} of $X$. And the map 

\begin{align*}
    \chi^m \mapsto x^{<m>}& = \prod_\rho x_\rho^{<m, u_\rho> } \end{align*}

induces the isomorphism

\begin{align*}
    \pi_\sigma^*: \C[\sigma^\vee \cap M] &\cong ( S_{x^{\hat{\sigma}}})^G \text{ for }G \text{ corresponding to the GIT quotient of } \mathbb{A}^{|\Sigma(1)|}_{ss} \\
        &\cong (S_{x^{\hat{\sigma}} } )_0 \end{align*}

which are degree 0 in the $Cl(X)$-grading on $S$. \\

So locally, we have 

\begin{align*}
(S_{x^{\hat{\sigma}}})_0 \otimes_k ( S_{x^{\hat{\sigma}} })_0 \stackrel{\psi}{\rightarrow} (S_{x^{\hat{\sigma}}})_0 \text{ by }\\
\sum \textbf{u}^\textbf{a} \otimes_k \textbf{u}^\textbf{b} \mapsto \sum \textbf{u}^{\textbf{a}+\textbf{b}}. \end{align*} for $\textbf{u}$ local coordinates in $U_\sigma$. \\

From 

\begin{align*}
I_X &= \left< \textbf{x}^\textbf{u} \otimes \textbf{x}^\textbf{v} - \textbf{x}^\textbf{v} \otimes \textbf{x}^\textbf{u} \text{ }|\text{ } \pi(\textbf{u})=\pi(\textbf{v}) \in Cl(X) \right> \\
&= \ker(\psi) \end{align*} for $\psi: R \otimes_k R \rightarrow k[Cl(X)]\otimes_k R \text{ via } \textbf{x}^\textbf{u}\otimes \textbf{x}^\textbf{v} \mapsto [\textbf{u}] \otimes \textbf{x}^{\textbf{u}+\textbf{v}}$, we have that locally:

\begin{align*}
(I_\Delta)_{x^{\hat{\sigma}} \otimes_k x^{\hat{\sigma}}} &= \ker( \textbf{x}^\textbf{u} \otimes_k \textbf{x}^\textbf{v} \mapsto [\textbf{u}]\otimes_k \textbf{x}^{\textbf{u}+\textbf{v}} ) \subset R_{x^{\hat{\sigma}}} \otimes_k R_{x^{\hat{\sigma}}} \end{align*}

Now 

\begin{align*}
    R_{x^{\hat{\sigma}}} \cong k[ x_1 , \dots, x_r, w_1^{\frac{+}{}}, \dots, w_\ell^{\frac{+}{}} ] \end{align*}
    
from getting powers of variables to cancel when localizing by  $x^{\hat{\sigma}}$ for some  $r$ and $\ell$, and

\[ 
(R_{x^{\hat{\sigma}}})_0 = 
k[\frac{x_1}{\textbf{w}^{\textbf{$\alpha$}_1}},
\dots,
\frac{x_r}{\textbf{w}^{\textbf{$\alpha$}_r}}
] \]

from $X$ unimodular implies that for each $b_k \in \sigma(1)$ there exists $m_k \in M$ such that $m_k \cdot b_i = \delta_i^k \text{ for all } b_k \in \sigma(1)$. Taking degree zero and considering the kernel gives 

\begin{align*}
(I_\Delta)_{x^{\hat{\sigma}} \otimes_k x^{\hat{\sigma}}} &= \left< \left\{ \frac{x_i}{w^{\alpha_i}} \otimes 1 - 1 \otimes \frac{x_i}{w^{\alpha_i}}  \right\}_{i=1}^r \right> \end{align*} which we'll write as

\begin{align*} \left< \{ x^u - y^u \} \right> \end{align*} on affine charts. \end{proof} 


Now\cite{bayer-popescu-sturmfels}, the ideal $I_X \subset S$ defining the diagonal embedding $X \subset X \times X$ equals the Lawrence ideal $J_L$ for the lattice $L=\text{Im}(B) = \ker(\pi)$ of principal divisors from the sequence at the beginning of this section~\ref{seq:toric.fund.seq.unimod}. This latice $L$ is unimodular precisely when $X_\Sigma$ is unimodular, and since the Lawrence lifting 

\[ \Lambda(L) = \{ (u, -u)\text{ }|\text{ }u\in L \} \]

is isomorphic to $L$, this implies that $\Lambda(L)$ is unimodular as well\cite{bayer-popescu-sturmfels}. \\

Next, we'll require the following definitions:\\

\textbf{Definition}: Given a field $k$, we consider the \underline{Laurent polynomial ring} $T = k[x_1^{\frac{+}{}}, \dots, x_n^{\frac{+}{}}]$ as a module over the polynomial ring $S = k[x_1, \dots, x_n]$. The module structure comes\cite{bayer-sturmfels} from the natural inclusion of semigroup algebras $S = k[\N^n] \subset k[\Z^n] = T$. \\

\textbf{Definition}: A \underline{monomial module} is an $S$-submodule of $T$ which is generated by monomials $\textbf{x}^\textbf{a} = x_1^{a_1}\cdots x_n^{a_n}, \textbf{a}\in \Z^n$. \\

\textbf{Definition}: When the unique generating set\cite{bayer-sturmfels} of monomials of a monomial module $M$ forms a group under multiplication, $M$ coincides with the \underline{lattice module} 

\[ M_L := S\{\textbf{x}^\textbf{a} \text{ }|\text{ }\textbf{a} \in L \} = k\{\textbf{x}^\textbf{b} \text{ }|\text{ } \textbf{b} \in \N^n + L \} \subset T \] 

for some lattice $L \subset \Z^n$ whose intersection with $\N^n$ is the origin $\textbf{0}=(0,0,\dots,0)$. \\

\textbf{Definition}: The $\Z^n/L$ graded \underline{lattice ideal} $I_L$ is\cite{miller-sturmfels}

\[ I_L = \left< \textbf{x}^\textbf{a} - \textbf{x}^\textbf{b} \text{ }|\text{ } \textbf{a}-\textbf{b} \in L \right> \subset S. \]

\textbf{Claim}: $I_{\Lambda(L)} = J_L$. \\

\begin{proof}

\begin{align*}
    I_{\Lambda(L)} &= \left< \mathbb{x}^\textbf{a}-\mathbb{x}^\textbf{b} \text{ }|\text{ }\textbf{a}-\textbf{b}\in \Lambda(L) \right>\\
    &= \left<\textbf{x}^{\textbf{a}_1}\textbf{y}^{\textbf{a}_2}-\textbf{x}^{\textbf{b}_1}\textbf{y}^{\textbf{b}_2}\text{ }|\text{ }\textbf{a}=(\textbf{a}_1,\textbf{a}_2), \textbf{b}=(\textbf{b}_1, \textbf{b}_2), \textbf{a}-\textbf{b} \in \Lambda(L) \right>\\
&= \left< \textbf{x}^{\textbf{a}_1}\textbf{y}^{\textbf{a}_2} - \textbf{x}^{\textbf{b}_1}\textbf{y}^{\textbf{b}_2}\text{ }|\text{ }(\textbf{a}_1-\textbf{b}_1, \textbf{a}_2-\textbf{b}_2) = (\textbf{a}_1,\textbf{a}_2) - (\textbf{b}_1,\textbf{b}_2) \in \{(\textbf{u},-\textbf{u})\text{ }|\text{ }\textbf{u}\in L \} \right>
\end{align*}

where the rightmost constraint implies $$ \textbf{a}_2-\textbf{b}_2 = -(\textbf{a}_1-\textbf{b}_1)=\textbf{b}_1-\textbf{a}_1 \text{ and } \textbf{a}_1-\textbf{b}_1\in L$$ gives

\begin{align*}
    &= \left<\textbf{x}^\textbf{u}\textbf{y}^\textbf{v} - \textbf{x}^\textbf{v} \textbf{y}^\textbf{u} \text{ }|\text{ } \textbf{u} - \textbf{v} \in L \right>= J_L
\end{align*}

\end{proof}

While the hull resolution of the unimodular Lawrence ideal $I_{\Lambda(L)}=J_L$ is not necessarily minimal, the minimal resolution of $I_{\Lambda(L)}$ does come from a cellular resolution of $M_{\Lambda(L)}$ and is described by a combinatorial construction given by the infinite hyperplane arrangement\cite{bayer-popescu-sturmfels} $\mathcal{H}_L$. Here, $\mathcal{H}_L$ comes from intersecting $\R L$, the real span of the lattice $L = \ker(\pi) \subset \Z^n$, with all lattice translates of the coordinate hyperplanes $\{x_i=j\}, 1\leq i \leq n, j\in \Z$. Bayer-Popescu-Sturmfels\cite{bayer-popescu-sturmfels} prove that each lattice point in $L$ is a vertex of the affine hyperplane arrangement $\mathcal{H}_L$, and that there are no additional vertices in $\mathcal{H}_L$ iff $L$ is unimodular. Furthermore, the quotient complex $\mathcal{H}_L/L$ of the hyperplane arrangement modulo $L$ supports\cite{bayer-popescu-sturmfels} the minimal $S$-free resolution of $J_L$, precisely because $\R L \cap \mathcal{H}_L = L$.

\subsection{Example of cellular resolution of $\mathcal{H}_L/L$ for $\PP^2$}

For the unimodular toric variety $X_\Sigma=\PP^2$, we have the fundamental toric exact sequence defining $Cl(X_\Sigma)$ via 

\[ 0 \rightarrow \Z^2 \stackrel{B}{\rightarrow} \Z^3 \stackrel{\pi}{\rightarrow} \Z \rightarrow 0 \]

with \[ B = \left[ \begin{matrix} 1 & 0 \\ 0 & 1 \\ -1 & -1 \end{matrix} \right] \]

and $B(a,b)$ given by 

\begin{align*}
\left[\begin{matrix} 1 & 0 \\ 0 & 1 \\ -1 & -1 \end{matrix}\right]\left[\begin{matrix} a \\ b \end{matrix}\right] &= \left[\begin{matrix} a \\ b \\ -a-b \end{matrix}\right]  \end{align*}

so that $\ker(\pi)=\text{Im}(B)$ contains

\[ \{e_1 - e_2, e_1-e_3, e_2-e_3 \} \] but can be generated by \[ \left\{ e_1-e_2, e_2-e_3 \right\}. \] 

Now $L = \Z \cdot \left< \left[\begin{matrix} 1 \\ -1 \\ 0 \end{matrix}\right], \left[\begin{matrix} 0 \\ 1 \\ -1 \end{matrix}\right] \right> $ is rank 2. \\

If we label $\mathcal{H}_L$ with elements of the monomial module $M_{\Lambda(L)}$ and build the corresponding quotient cellular complex on $\mathcal{H}_L/L$ labeled with one up triangle and one choice of down triangle, subject to the orientation below, we see\cite{bayer-popescu-sturmfels}: \\[1cm]

\begin{figure}[h]
\includegraphics[width=\textwidth]{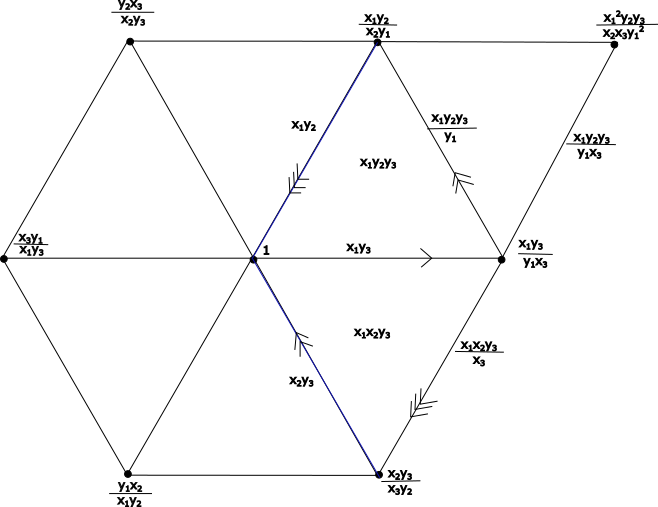}
\caption{Cellular Complex $\mathcal{H}_L/L$ for $\PP^2$}
\label{fig:bpsP2Fig1}
\end{figure}

Here, the edges corresponding to the basis vectors $e_1-e_2$ and $e_2-e_3$ are labeled in blue, the edge ordering is given by the number of arrows (1-3), and the upward-pointing face of $\mathcal{H}_L/L$ is oriented counterclockwise while the downward-facing face is oriented clockwise in Figure~\ref{fig:bpsP2Fig1}. We order the 2-cells via $\{\text{``up-triangle," ``down-triangle"}\}$. The cellular complex $\mathcal{H}_L/L$ has one vertex, three edges, and two 2-cells. This choice of fundamental domain and face ordering gives the S-free resolution\cite{bayer-popescu-sturmfels}  $(\textbf{F}_{\mathcal{H}_L/L}^\bullet, \partial)$ resolution of $S/J_L$:\\

\begin{equation}\label{eqn:3.2}
\begin{array}{ccccccccc}
0 &\rightarrow & S^2 &\stackrel{\left(\begin{matrix} y_2 & x_2 \\ y_1 & x_1 \\ y_3 & x_3 \end{matrix}\right) }{\longrightarrow}& S^3 &\stackrel{\left(\begin{matrix} y_1x_3 - x_1y_3 \\ x_2y_3 - x_3y_2\\ x_1y_2 - x_2y_1 \end{matrix}\right) }{\longrightarrow}& S^1 & \rightarrow & 0 \\ & & & & & & \downarrow\\ & & & & & & \faktor{S}{J_L} & \rightarrow & 0 
\end{array} 
\end{equation}

where

\[\partial(F) = \sum_{\makecell{F'\subset F,\\ codim(F',F)=1}}\epsilon(F,F')\frac{m_F}{m_{F'}} \cdot F' \]

for $m_F$ the monomial label of the face $F$ and $\epsilon(F,F')\in \{\frac{+}{}1\}$ with $\epsilon(F,F')=1$ if the orientation of $F'$ agrees with that induced by $F$. Here, $m_F$ is the lowest common multiple of $m_{F'}$ over the facets of $F$, starting with a labeling of the vertices. \\

That is, the cellular resolution of $\faktor{S}{J_L}$ coming from $M_{\Lambda(L)}$ in Equation~\ref{eqn:3.2} resolves the ideal $I_{\PP^2}$ of the diagonal $\PP^2 \subset \PP^2\times \PP^2$. \\

\subsection{Discussion}
Here, the cellular resolution for $M_{\Lambda(L)}$ works precisely because the vertices of a unimodular lattice are exactly the vertices of $\mathcal{H}_L$. In order to resolve the diagonal for a toric D-M stack associated to a given simplicial toric variety, we require the Morita theory of Section~\ref{subsec.3.4.6}.

We generalize the Cellular Resolution of $\mathcal{H}_L/L$ for toric D-M stacks associated to some simplicial toric varieties in Section~\ref{subsec:AdaptArgSimplicial}.




\section{Resolution of diagonal for global quotient of unimodular toric variety by finite abelian group} 
\subsection{Resolution of diagonal in bimodules}

Let $A$ be an algebra over $R$, a commutative algebra over $k$ a field. $A$ is not necessarily commutative. We consider $I_\Delta \leq A \otimes_k A$ as a submodule given by $\ker \mu: A\otimes_k A \rightarrow A$ of multiplication. Here, $A \otimes_k A \geq A$ contains an isomorphic copy of $A$. We consider $A \otimes_k A$ as an $A \otimes_k A^{op}$-bimodule, not as an algebra. Here, $$ A_\Delta = \faktor{A \otimes_k A}{I_\Delta}$$ 
    
    To consider derivations and differentials of elements in $A$, consider the map given by
    
    \[ \alpha \in A \mapsto \tilde{d}\alpha = \alpha \otimes 1 - 1 \otimes \alpha \in I_\Delta \]
    
    We do not consider $I/I^2$ here. \\
  
  Claim: $\ker\mu = I_\Delta \cong <\tilde{d}\alpha \text{ }|\text{ }\alpha \in A>$ as $A\otimes A^{op}$-module.
  
  \begin{proof}
  ($\Rightarrow$): We construct $g: \ker \mu \rightarrow \left< \tilde{d}\alpha: \alpha\in A\right>$ via $g(\sum_i z_i\otimes_k w_i) = \sum_i z_i \tilde{d}w_i$ as follows:\\
  
  \begin{align*}
      \sum_i z_i\otimes w_i &\in \ker \mu \implies \sum_i z_iw_i = 0 \text{ so that }\\
      \sum_i z_i\otimes w_i &= \sum_i z_i(1\otimes w_i - w_i\otimes 1) \in \left< \tilde{d}\alpha \text{ }|\text{ }\alpha\in A \right>. \end{align*}

($\Leftarrow$): We clearly have that $\tilde{d}\alpha = 1\otimes_k \alpha - \alpha\otimes_k 1 \in \ker\mu$. 
  \end{proof}

    Relative to $R$, if we suppose that $A$ is finitely generated as an $R-$module over $R$ by $\{f_1, \dots, f_r\}$ then \[ \left< \tilde{d}a \text{ } |\text{ }a \in A \right> \cong \left< \tilde{d}r \text{ }|\text{ }r\in R \right> \otimes \{f_1, \dots, f_r \} \]
    so that this gives a finite set of generators for $I_\Delta$. \\

Additionally, we observe that the Leibniz rule     $\tilde{d}(\alpha\beta) = \tilde{d}(\alpha)\beta + \alpha\tilde{d}(\beta)$ holds here.

\subsection{Constructing local diagonal object via Morita theory}
\label{subsec.3.4.6}
Consider $\faktor{\A^n}{\mu}$ where $\mu$ is a finite group of the torus. We have that $D(\faktor{\A^n}{\mu}) = D(k[z_1,\dots,z_n]-\text{Mod}^\mu)$ where modules are considered with grading. The goal here is to understand $D(\faktor{\A^n}{\mu} \times \faktor{\A^n}{\mu})$ in terms of $D(\faktor{\A^n}{\mu})$, and more specifically, to understand $\mathcal{O}_\Delta$ of $D(\faktor{\A^n}{\mu}\times\faktor{\A^n}{\mu})$ in terms of $\mathcal{O}_\Delta$ on $\faktor{\A^n}{\mu}$. That is, we would like a description of $\mathcal{O}_\Delta$ of $\faktor{\A^n}{\mu}\times \faktor{\A^n}{\mu}$ in terms of $\Delta$ for affine space. There is always the cohomological grading present, but now we also have an extra grading on $k[\A^n]-\text{Mod}^\mu$. Here, let $S = k[\A^n]$. Note that $D(\faktor{\A^n}{\mu})$ is generated by line bundles $\mathcal{O}(\xi)$ for characters $\chi\in \mu^*$, so we have the generator

\[ G = \sum_{\chi \in \mu^*}\mathcal{O}(\chi) \] 
    
Note: $Hom(\mathcal{O}(\chi_i),\mathcal{O}(\chi_j))$ is $S$ with the right grading. Ignoring grading, we have 

\begin{align*}
    G &= \sum_{|\mu|}  S \end{align*}  so that 
\begin{align*}    End^\bullet(G) &= Mat_{|\mu|}(S).
\end{align*}
    
    In order to incorporate the grading, we would then take the degree $0$ part. Here we use Morita theory.

\subsection{Background on Morita Equivalence}

   Note: Notation in this sub-subsection follows Weibel, and differs from notation used in later sections.\\
    \textbf{Definition}: Two unital rings $R,S$ are said to be \textbf{Morita equivalent} if there exists an $R-S$ bimodule $P$ and $S-R$ bimodule $Q$ such that $P\otimes_S Q \cong R$ as $R-R$ bimodules and $Q\otimes_R P \cong S$ as $S-S$ bimodules. Then \\
    \begin{center}
    $(-)\otimes_R P$ gives a functor from $\text{mod}-R \rightarrow \text{mod}-S$ and\\
    $(-)\otimes_S Q$ gives a functor from $\text{mod}-S \rightarrow \text{mod}-R$\\
    \end{center}
    
    which are inverse categorical equivalences, since for all $M\in \text{ mod }-R$, we have
    
    \begin{align*}
        M\otimes_R P \otimes_S Q &\cong M\otimes_R (P\otimes_S Q) \cong M \end{align*}
        
        and
        
        \begin{align*}
        M\otimes_S Q \otimes_R P &\cong M\otimes_S (Q\otimes_R P)\cong M
    \end{align*}
for all $M\in \text{mod}-S$.\\


\subsection{Choosing cofinite lattice $\tilde{L}\subset L$ for global quotient of unimodular toric variety}
\label{subsec:AdaptArgSimplicial}

For the toric D-M stack associated to $X_\Sigma'$, a simplicial toric variety which is a global quotient $\faktor{X_\Sigma}{G}$ of a unimodular toric variety $X_\Sigma$ with $G$ a finite abelian group, we can choose an isomorphism of the lattice $L$ of principal divisors of $X_\Sigma$ which is unimodular such that $\tilde{L}\subset L$ is cofinite, and $\faktor{L}{\tilde{L}} \cong G$. Furthermore, locally, $X_\Sigma' \cong \faktor{\A^n}{G}$, so that locally we can consider $L \cong \Z^n$.  

\subsection{Diagonal Object in Mod-A}
\begin{figure}[h]
\centering
\includegraphics[width=10cm]{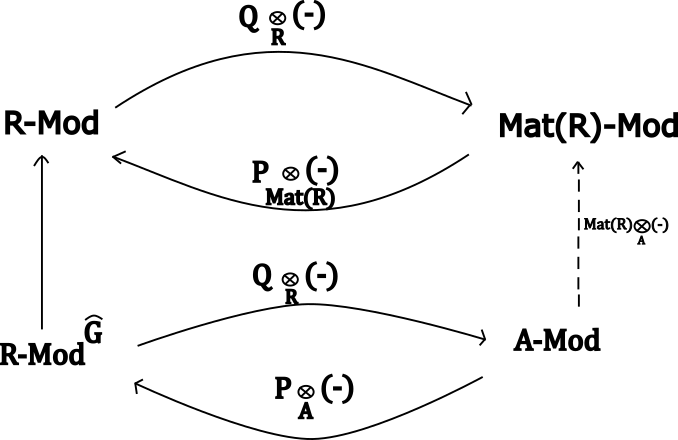}
\caption{Morita Equivalence Functors Between $R-Mod$ and $Mat(R)-Mod$, and their Graded and $\hat{G}$-Invariant Counterparts Between $R-Mod^{\hat{G}}$ and $A-Mod$}
\label{fig:3.5.7}
\end{figure}

Here, let $R = k[\A^n]$ and $S = R\otimes_k R \cong k[x_1,\dots,x_n,y_1,\dots,y_n]$. 
In $D(\faktor{\A^n}{G})$, we have generator $H=\bigoplus_{\chi\in\hat{G}}\mathcal{O}(\chi)$, with endomorphism algebra $A=End(H)$. Now, $A\simeq \Delta \in A\otimes A^{op}-Mod$. We have functors given in Figure~\ref{fig:3.5.7}

\flushleft 

Here, $P \in (R^{\hat{G}}, Mat(R))-Mod$ is a row vector with grading on generators 

\[ P = \left< e_\chi\text{ }|\text{ } |e_\chi|=-\chi\right>, \] 

$Q\in (Mat-R, R^{\hat{G}})-Mod$ is a column vector with grading on generators

\[ Q = \left< e_\chi\text{ }|\text{ }|e_\chi|=\chi\right> \] 

and we have 

\[ Hom_{Mat}(Q,-) \simeq P\otimes_{R^{\hat{G}}} ( - ) \] 

For $P\otimes_A (-): A-Mod \rightarrow R-Mod^{\hat{G}}$ applied to $\Delta \cong A \in Obj(A\otimes A^{op}-Mod)$, we have that 

\[ P\otimes_A A = \left< z^\alpha e_{\rho_1\rho_2} \text{ }|\text{ } |z^\alpha e_{\rho_1\rho_2}| = |z^\alpha| + \rho_1 - \rho_1 + \rho_2 \right> \] 

with grading

\begin{align*}
|z^\alpha e_{\rho_1\rho_2}| &= |z^\alpha| + \rho_1-\rho_1 + \rho_2\\
&= |z^\alpha| + \rho_2. \end{align*}

For $\pi(M_{\Lambda(L)}) = M_{\Lambda(L)} \otimes_{S[\Lambda(\tilde{L})]} S$, we have 

\begin{align*}
\pi(M_{\Lambda(L)}) &= S \cdot \{ \frac{x^u}{y^u} \text{ }|\text{ } u \in L \} \otimes_{S[\Lambda(\tilde{L})]} S
\end{align*}

with $S[\Lambda(\tilde{L})]= \left< x^{u_1}y^{u_2} \otimes z^{(v,-v)} \text{ }|\text{ } (u_1,u_2) \in \mathbb{N}^{2n}, v\in \tilde{L} \right>$ induces an $\faktor{L}{\tilde{L}}$ grading left on, say, $\{x_1,\dots,x_n\}$ from cancelling $y^u$ in the denominator. Here, the map

\[ \phi: S[L]\rightarrow M_L \text{ via } x^{u_1}y^{u_2} \otimes z^{(v_1,v_2)} \mapsto x^{u_1 + v_1}y^{u_2+v_2} \] induces an $S[\Lambda(L)]$-module structure on $M_L$, hence an $S[\Lambda(\tilde{L})]$-structure on $M_{\Lambda(L)}$. In Bayer-Sturmfels\cite{bayer-sturmfels}, the map $\phi: S[L]\rightarrow M_L$ was used to show that 

\[ \pi(M_L) \cong \faktor{S[L]}{\ker\phi} \otimes_{S[L]} S \cong \faktor{S}{I_L}\] 

by using the fact that $\pi(M_L)$ is cyclic when we only have the unimodular lattice $L$ to consider. \\[1cm]

We can also view

\[ (R-Mod)^{\hat{G}} \rightarrow Mod-A \] by

\begin{align*}
\bigoplus_{\rho\in\hat{G}} M_\rho \cong M_\cdot \mapsto Hom(\bigoplus_{\rho\in \hat{G}} R(-\rho), M) \cong \bigoplus_\rho Hom(R(-\rho),M) \cong \bigoplus_\rho M_\rho \end{align*}

where $\bigoplus_\rho M_\rho$ contains idempotents $e_{\rho\rho}\in A$ on the right-hand side. 

We also have $\pi(M_{\Lambda(L)}) \in (S-Mod)^{\hat{G}\times\hat{G}}$ with

\begin{align*}
\pi(M_{\Lambda(L)}) \simeq +_{\overline{u}\in\hat{G}} M_{\tilde{L}} \cdot \tilde{z}^{\overline{u}} 
\end{align*}

is a sum which is not necessarily direct, and 

\begin{align*}
S\cdot \left\{ \frac{x^u}{y^u} \text{ }|\text{ }u\in L \right\} \subset S_{(x_1\dots x_n y_1 \cdots y_n )} \cong k[\Z^{2n}]. \end{align*}

Here, the localized ring $S_{(x_1\dots x_n y_1\dots y_n)}$ admits a $\Z^{2n}$ grading descending to an $\faktor{L}{\tilde{L}} \oplus \faktor{L}{\tilde{L}}$ grading. While $\pi(M_{\Lambda(L)})$ is no longer cyclic, we do have

\begin{align*}
S[L] \cong \bigoplus_{\overline{u} \in \hat{G}} S[\tilde{L}]\cdot z^{\overline{u}} \end{align*}

and

\begin{align*}
S[\Lambda(L)] \cong \bigoplus_{u\in \hat{G}} S[\Lambda(\tilde{L})] \cdot z^{(u,-u)} \end{align*}

carries an anti-diagonal grading. Now $$S[L] = \left< x^u \otimes z^v \text{ }|\text{ } v \in L\right>$$ and in $$ S[\Lambda(L)] = \left< x^{u_1}y^{u_2} \otimes z_1^u z_2^{-u} \text{ }|\text{ } u \in L \right>$$ $S$ has an $L \times L$ grading since $L\cong \Z^n$, and taking the quotient of $L$ by $\tilde{L}$ induces a quotient on $\Lambda(L)$ by $\Lambda(\tilde{L})$ which gives an $\faktor{L}{\tilde{L}} \times \faktor{L}{\tilde{L}}$-grading on $S[L]$. Now $S[\tilde{L}]$ has a degree 0 grading on $z'$s with respect to $\faktor{L}{\tilde{L}}$. \\[.2cm]

For $\phi: S[L]\rightarrow M_L$ via $x^u \otimes z^v \mapsto x^{u+v}$, we have that $\phi$ is surjective, $\ker\phi = \left< x^u \otimes z^v - x^{u+v}\otimes 1\right>$ and $$I_L = \left< x^u - x^v \text{ }|\text{ }u-v\in L \right> \leq S $$ comes from $\ker\phi$ by setting $z's$ equal to 1. Here, we have that 

\begin{align*}
M_L \cong \faktor{S[L]}{\ker\phi} \end{align*} as $S[L]$ modules. As $S[\tilde{L}]$-modules, we still have 

\[ S[L] \cong \bigoplus_{\overline{u}\in\hat{G}} S[\tilde{L}] \cdot z^{\overline{u}} \] 

For $\phi: S[\Lambda(L)] \rightarrow M_{\Lambda(L)}$, we have that 

$$ \ker\phi = \left< y_i \otimes z^{(e_i, -e_i)} - x_i \right>$$ is homogeneous of degree $(1,0)$ is one generating set, and $\ker\phi = \left< x_i \otimes z^{(-e_i, e_i)} - y_i\right>$ is another generating set of homogeneous elements of degree $(0,1) \in (\hat{G}, \hat{G})$. Furthermore, we have that 

\begin{align*}I_{\Lambda(L)} = \left< x_j - y_j \right> \end{align*} since $L \cong \Z^n$, and considering $\phi$ as a map of $S[\tilde{L}]$-modules gives that $\ker\phi$ still specializes to $I_{\Lambda(L)}$ when we set the $z$'s equal to $1$. \\[2cm]

To show that 

\[ \begin{array}{ccccc}
A \cong &\pi(M_{\Lambda(L)}) = M_{\Lambda(L)} \otimes_{S[\Lambda(\tilde{L})]} S \in (S-Mod)^{\hat{G}\times\hat{G}}\\
& || \\
& M_L \otimes_{S[\tilde{L}]} S \end{array}
\]

we note that $A$ is generated over $k$ by 

\[ A/k = \left< w^\alpha e_{\rho_1\rho_2} \text{ }|\text{ } |\alpha| = \rho_1 - \rho_2 \right> \]

so we construct a map $\Psi: A \rightarrow \pi(M_{\Lambda(L)})$ by describing its image on generators of $A/k$:\\

\begin{align*}
\psi(w^\alpha e_{\rho_1 \rho_2}) = x^\alpha \otimes z^{\Lambda(\rho_2)}. \end{align*}

Given $M\in R-Mod$, we can use the projectors $e_{\rho_i\rho_i}$ to get \[ e_{\rho_i\rho_i} M = M_{\rho_i} \in k-v.s. \] and

\[ M = \bigoplus_i M_{\rho_i} \] as a graded $R-$module.

To get a $\hat{G}\times \hat{G}$-grading on $A\in (A\otimes A^{op})-Mod$, we use projectors \[ e_{\rho \rho} A e_{\rho'\rho'}\] to pick up the $(\rho, -\rho')$-homogeneous component of $A$ in the $\hat{G}\times \hat{G}$-grading. \\

To get a $\hat{G}\times \hat{G}$-grading on $\pi(M_L)$, we have that 

\[ |x^{\alpha_1}y^{\alpha_2}\otimes z_1^\beta z_2^{-\beta}| = (\alpha_1+\beta, \alpha_2-\beta). \]

\begin{figure}[h]
\centering
\includegraphics[width=.5\textwidth]{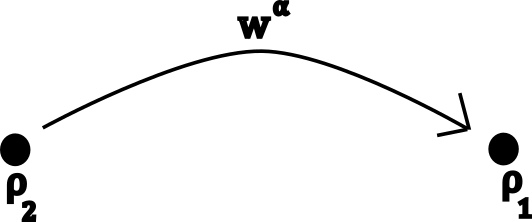}
    \caption{Quiver associated to $w^\alpha e_{\rho_1\rho_2}$}
    \label{fig:my_label}
\end{figure}

Here, the grading on generators of $A$ over $k$ comes from Figure~\ref{fig:my_label}.\\[1cm]

Since $\pi(M_{\Lambda(L)}) \in (S-Mod)^{\hat{G}\times\hat{G}}$ and $A \simeq \Delta \in (A \otimes A^{op})-Mod$, we put an $(S-Mod)^{\hat{G}\times \hat{G}}$ structure on $A$ via

\[ S = R\otimes_k R \] 

and use the matrix representation representing multiplication by $r\in R$ to give $R\hookrightarrow A$ so that both $A$ and $\pi(M_{\Lambda(L)})\in (S-Mod)^{\hat{G}\times \hat{G}}$. 

Now to check that $\psi$ is compatible with the $(S-\text{ Bimod })^{\hat{G}\times\hat{G}}$ structure, we note that 

\begin{align*}
x^\beta w^\alpha e_{\rho_1\rho_2} &= w^{\alpha+\beta} e_{\rho_1+|\beta|, \rho_2} & \text{ Left action}\\
\end{align*}

and

\begin{align*}
y^\beta w^\alpha e_{\rho_1\rho_2} &= w^{\alpha+\beta} e_{\rho_1, \rho_2-|\beta|} &\text{ Right Action} 
\end{align*}

gives that we have an $S-Mod$ homomorphism: \[ f(s\cdot m) = s\cdot f(m) \] 

since 

\begin{align*}
\psi(x^\beta w^\alpha e_{\rho_1\rho_2}) &= \psi(w^{\alpha+\beta}e_{\rho_1+|\beta|, \rho_2} )\\
&= x^{\alpha+\beta} \otimes z^{\lambda(\rho_2)}\\
&= x^\beta \cdot (x^\alpha \otimes z^{\lambda(\rho_2)}\\
&= x^\beta \psi(w^\alpha e_{\rho_1\rho_2} )\\
\end{align*}

and

\begin{align*}
\psi(y^\beta w^\alpha e_{\rho_1\rho_2}) &= \psi(w^{\alpha+\beta} e_{\rho_1, \rho_2-|\beta|})\\
&= x^{\alpha+\beta}\otimes z^{(\rho_2-|\beta|, -\rho_2+|\beta|)}\\
&= (x^\beta\otimes z^{(-\beta,\beta)})\cdot(x^\alpha \otimes z^{\lambda(\rho_2)})\\
&= y^\beta (x^\alpha\otimes z^{\lambda(\rho_2)})\\
&= y^\beta \psi(w^\alpha e_{\rho_1\rho_2} ) \end{align*} 

Lastly, to show that this is an isomorphism, we construct an inverse to $\Psi$. Since $\pi(M_{\Lambda(L)}) \cong R[\hat{G}]$ as $k-v.s.$, we construct

 \[ \Phi: x^\alpha \otimes z^{\overline{v}} \mapsto w^\alpha e_{\alpha+\overline{v},\overline{v}} \] 

as an inverse to $\psi$. \\

This shows that for global quotients, we tensor over $S[\Lambda(\tilde{L})]$ rather than over $S[\Lambda(L)]$ to get $\Delta$.

\subsection{Global diagonal object and resolution of the diagonal for $X_\Sigma' = \faktor{X_\Sigma}{G}$ a global quotient of unimodular $X_\Sigma$ by $G$ a finite abelian group}

Here, we'll give a global diagonal object for the example of $\faktor{\PP^n}{G}$ and a proof for why this is a global diagonal object for $\faktor{\PP^n}{G}$ which depends only on the fact that $\PP^n$ is unimodular, hence giving a global diagonal object for the more general case of $X_\Sigma' = \faktor{X_\Sigma}{G}$ a global quotient of unimodular $X_\Sigma$ by $G$ a finite abelian group. \\[.2cm]

\textbf{Claim}: We have a global diagonal object for $\faktor{\PP^n}{G}$ for $G$ a finite abelian group. 

\begin{proof} Locally, $\faktor{\PP^n}{G}$ is isomorphic to $\faktor{\A^n}{G}$. So locally, we have $U_\sigma' = \faktor{U_\sigma}{G}$ for $U_\sigma = \A^n$ open in $X_\Sigma = \PP^n$. For $\A^n$, we have presentation of $Cl(\A^n)=0$ via 

\[ 0 \rightarrow M \stackrel{B}{\rightarrow} \Z^{\Sigma(1)} \stackrel{\pi}{\rightarrow} Cl(X) \rightarrow 0 \] 

gives 

\[ 0 \rightarrow \Z^n \stackrel{Id}{\rightarrow} \Z^n \stackrel{\pi}{\rightarrow} 0 \rightarrow 0 \] 

and $L= \ker(\pi)=Im(B)\cong \Z^n$. This gives $\Lambda(L)= \{(u,-u) \text{ }|\text{ } u\in \Z^n \}$. We have $\Delta$ for $\faktor{\A^n}{G}$ as

\begin{align*}
\Delta \cong \pi(M_{\Lambda(L)}) &= \faktor{S[\Lambda(L)]}{\ker\phi} \otimes_{S[\Lambda(\tilde{L})]} S \\
&= \faktor{\left( k[\A^{2n}]\otimes \left< z^{\Lambda(e_i)\frac{+}{}} \right> \right) }{\left< y_i - x_i \otimes z^{\Lambda(e_i)} \right> } \otimes_{S[\Lambda(\tilde{L})]} k[\A^{2n}]\\
&\cong R[\hat{G}]
\end{align*}

where $R[\hat{G}]$ is the group-ring with variables $w_i$ and $G \cong \faktor{L}{\tilde{L}}$. Now the transition maps for coordinates $x_i$ on local charts of $\A^n = U_\sigma \subset X_\Sigma=\PP^n$ induce transition maps on transition maps for $U_\sigma \times U_\sigma$, and hence on $U_\sigma' \times U_\sigma' = \faktor{U_\sigma}{G} \times \faktor{U_\sigma}{G}$. This induces transition maps on the $z^{\Lambda(e_i)}$'s from locally writing $x_i = y_i \otimes z^{\Lambda(e_i)}$ in $M_L$ in $U_\sigma'$. This implies that the transition maps for the $z^{\Lambda(e_i)}$'s glue together, and the diagonal object $\Delta$ for $\faktor{\A^n}{G}$ glues together to give a global diagonal object for $\faktor{\PP^n}{G}$. \end{proof} 

Note: This proof is highly dependent on the fact that locally, $\PP^n$ is isomorphic to $\A^n$ in the sense that $\PP^n = \cup_i U_{\sigma_i}$ with each $U_{\sigma_i}\cong \A^n$. Since this also holds for any unimodular toric variety by the \textbf{Lemma} below, an analogous argument gives a global diagonal object for $X_\Sigma' = \faktor{X_\Sigma}{G}$ with $X_\Sigma$ unimodular. \\[2cm]

\textbf{Lemma}: A complete unimodular toric variety $X_\Sigma$ has a cover by open charts which are all isomorphic to $\A^n$. \\

\begin{proof} Let the complete unimodular toric variety $X_\Sigma$ have fan $\Sigma$ with maximal-dimensional cone $\sigma \in \Sigma(n)$. Then $\sigma^\vee \cap M$ is generated as a semigroup by normal vectors to the facets of $\sigma$, which is to say by the ray generators of $\sigma^\vee$. But there are $n$ such normal vectors, and

$X_\Sigma$ unimodular impies that $\sigma $ is unimodular, so subsets of $\rho_i \in \sigma(1)$ up to size $n$ are linearly independent over $\Z$. This implies that the $n$ ray generators of $\sigma^\vee$ are also linearly independent over $\Z$. So $U_\sigma \cong \A^n$. \end{proof}

\textbf{Claim}: $\pi(M_{\Lambda(L)})$ is resolved by the cellular complex associated to $\faktor{\mathcal{H}_L}{\tilde{L}}$, which we view as \[ \faktor{\mathcal{H}_L}{\tilde{L}} \cong \bigoplus_{u\in G} \faktor{\mathcal{H}_L}{L}\]

together with a covering map of 

\[ \begin{array}{c}
\faktor{\mathcal{H}_L}{\tilde{L}}\\
\downarrow \\
\faktor{\mathcal{H}_L}{L} 
\end{array}\]

\begin{proof}
Since $M_{\Lambda(L)}$ is a monomial module, $\mathcal{H}_L$ resolves $M_{\Lambda(L)}$. The complex $(\mathcal{F}_{\mathcal{H}_L}^\bullet, \partial)$ is not $S$-finite, but has finite length $m=\text{rank}(L)$. It is also a minimal $\Z^{2n}$-graded free $S$-resolution of the lattice module $M_{\Lambda(L)}$ \cite{bayer-popescu-sturmfels}. 
Since $\faktor{\mathcal{H}_L}{L}$ resolves \[ M_{\Lambda(L)} \otimes_{S[\Lambda(L)]} S\]
we have that $\faktor{\mathcal{H}_L}{\tilde{L}}$ resolves \[ \pi(M_{\Lambda(L)}) = M_L \otimes_{S[\Lambda(\tilde{L})]} S. \]
\end{proof}







\subsection{Example of resolution for $\faktor{\PP^1}{\mu_6}$}

Recall that for $\PP^1$ we have a presentation of $Cl(\PP^1)\cong \Z$ given by 

\[\begin{array}{ccccccccc}
0 &\rightarrow & M & \stackrel{B}{\rightarrow} & \Z^{|\Sigma(1)|} &  \stackrel{\pi}{\rightarrow} & Cl(\PP^1) & \rightarrow 0\\
 & & ||\wr & & ||\wr & & ||\wr \\
 0 & \rightarrow & \Z & \rightarrow  & \Z^2 & \rightarrow & \Z & \rightarrow & 0
 \end{array}\] 

 with maps $B$ and $\pi$ given by

 \[ 0 \rightarrow \Z  \stackrel{ \left(\begin{matrix} 1 \\ -1 \end{matrix} \right)}{\rightarrow} \Z^2 \stackrel{(1 \text{ }1 ) }{\rightarrow} \Z \rightarrow 0 \] 

 so that $L = \ker(1 1)$ with $L = \Z\left< \left( \begin{matrix} 1 \\ -1 \end{matrix} \right) \right> \subset \Z^2$. If we fix $L \cong \Z$, then $\tilde{L} = 6\Z$ and $\faktor{L}{\tilde{L}}\cong \mu_6=G$. Now $\PP^1$ has monomial labels on $\faktor{\mathcal{H}_L}{L}$ given by Figure~\ref{fig:my_label2}:

\begin{figure}[h]
\centering
\includegraphics[width=.5\textwidth]{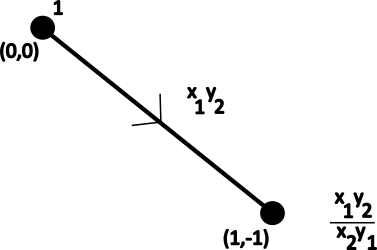}
    \caption{$\faktor{\mathcal{H}_L}{L}$ for $\PP^1$}
    \label{fig:my_label2}
\end{figure}

\begin{figure}[h]
\centering
\includegraphics[width=\textwidth]{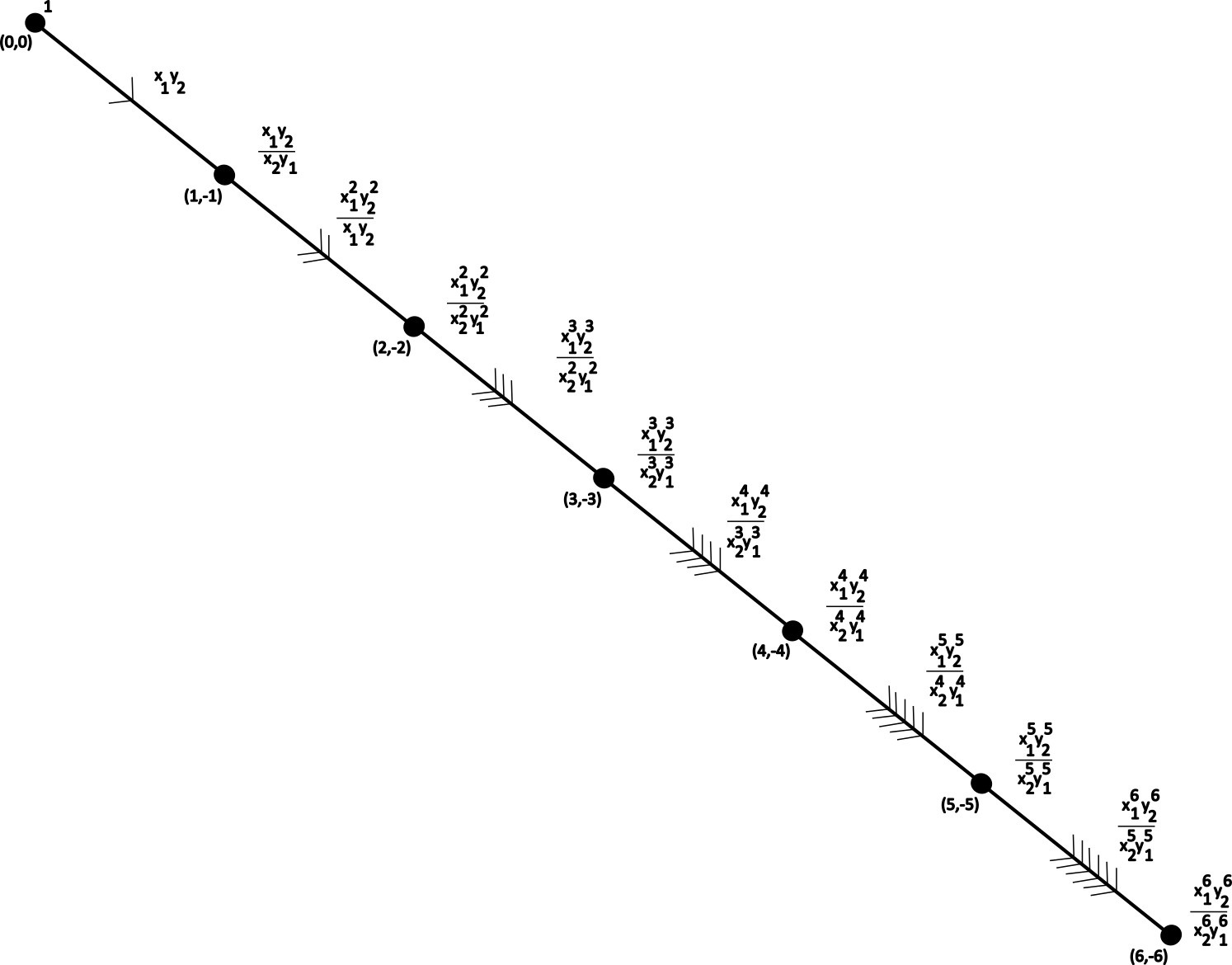}
    \caption{$\faktor{\mathcal{H}_L}{\tilde{L}}$ for $\faktor{\PP^1}{\mu_6}$}
    \label{fig:my_label3}
\end{figure}

 For $G = \mu_6$, we have the cellular complex $\faktor{\mathcal{H}_L}{\tilde{L}}$ and is given by Figure~\ref{fig:my_label3}.

\newpage

\begin{figure}[h]
\centering
\includegraphics[trim={4.5cm 5cm 5cm 4.5cm},clip, width=.5\textwidth]{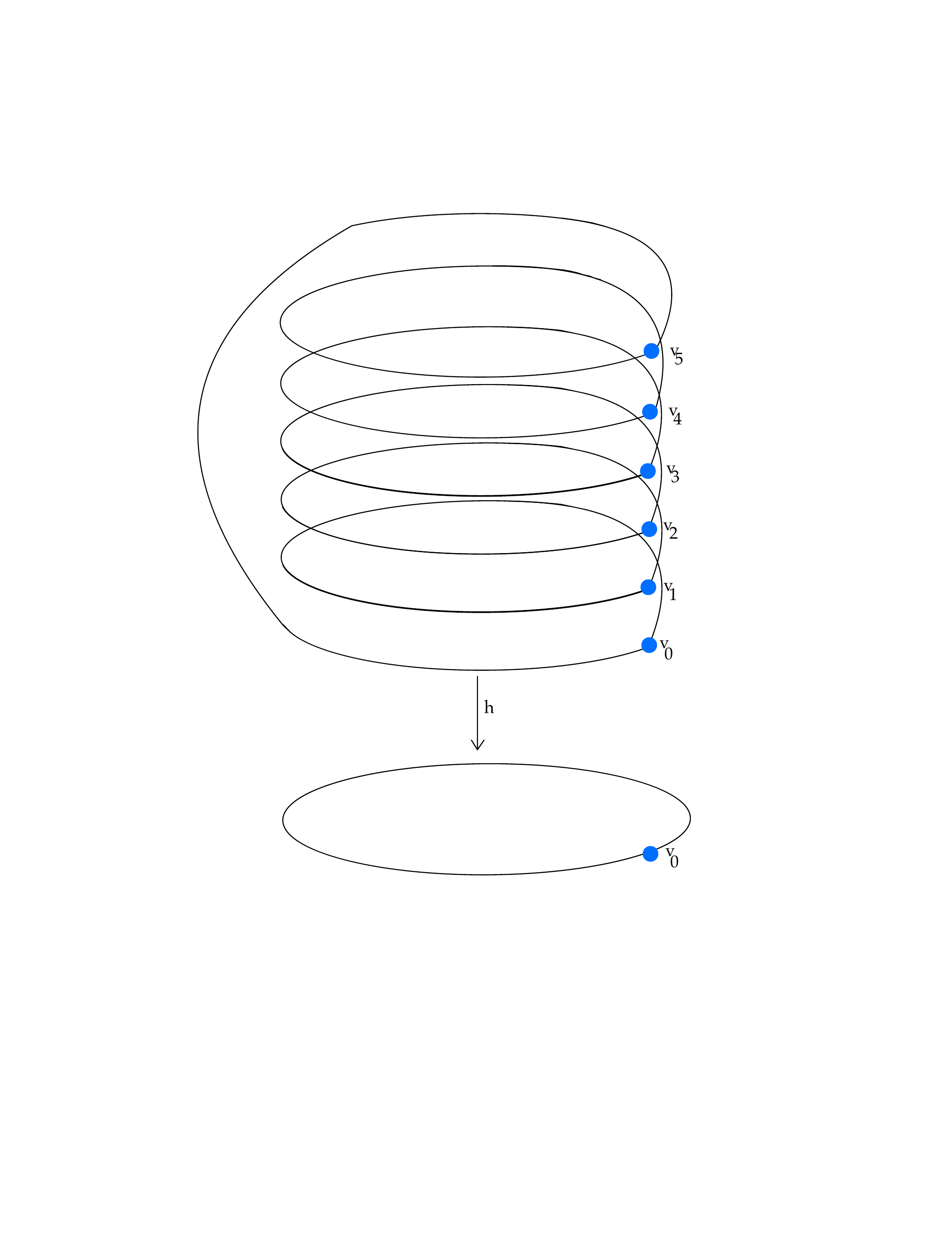}
    \caption{$\faktor{\mathcal{H}_L}{\tilde{L}}$ for $\faktor{\PP^1}{\mu_6}$}
    \label{fig:my_label4}
\end{figure}

The decomposition

\[ \faktor{\mathcal{H}_L}{\tilde{L}} = \bigoplus_{u\in G} \faktor{\mathcal{H}_L}{L} \] and corresponding covering map $h$ is given by Figure~\ref{fig:my_label4}
with corresponding complex and resolution 

\[
\begin{array}{ccccc}
S^6 & \stackrel{\left( \begin{matrix} -x_1y_2 & 0 & \cdots & x_2y_1 \\ x_2y_1 & -x_1y_2 & \cdots & 0 \\
0 & x_2y_1 & \cdots & 0 \\ 
\vdots & \vdots & \ddots & \vdots\\
0 & 0 & \cdots & -x_1y_2 \end{matrix} \right)}{\longrightarrow} & S^6 & \rightarrow & 0 \\
 & &  \downarrow & \\
 & &  \pi(M_{\Lambda(L)}) & \rightarrow & 0 \end{array} \]

The covering map $h$ corresponds to a subdivision of the maximal torus of $\PP^n$ by $G$.

\subsection{Upgrading the grading on $(S-Mod)^{\hat{G} \times \hat{G}}$}

If we add in a $(\Z\oplus \hat{G}) \times (\Z\oplus \hat{G})$-grading to $(S-Mod)^{\hat{G}\times \hat{G}}$, and view $\pi(M_{\Lambda(L)}) \in Ob\text{ }(S-Mod)^{\hat{G}\times \hat{G}}$, for instance, for $\faktor{\PP^1}{\mu_4}$ with $R = \C[x_0,x_1]$, exact sequence

\begin{align*}
0 \rightarrow \Z \stackrel{\left( \begin{matrix} 1 \\ -1 \end{matrix} \right)}{\rightarrow} \Z^2 \stackrel{\left( \begin{matrix} 4& 4 \end{matrix} \right)}{\rightarrow} \Z \rightarrow G \rightarrow 0 
\end{align*}

with map

\begin{align*}
\Z \oplus \faktor{\Z}{4\Z} \leftarrow \Z^2\\
(a+b,\overline{a}) \leftarrow (a,b) \end{align*}

and

\begin{align*}
|x_0^{a_0}x_1^{a_1}| = (a_0+a_1, \overline{a_0})
\end{align*}

with $a_0=-a_1$, $a_0=4b$ from $L=Im(B)\subset \Z^2$ gives a $\Z\oplus G$-grading on $R$, inducing a $(\Z\oplus G)\times(\Z\oplus G)$-grading on $S$. \\

Now, with \[ S[\Lambda(L)] \cong \bigoplus_{\overline{u} \in \hat{G}} S[\Lambda(\tilde{L})]\cdot z^{\Lambda(\overline{u})} \] 

allows us to write the resolution for $\faktor{\PP^1}{\mu^4}$ as

\[
\begin{array}{cccccc}
0 \rightarrow S^4 & \stackrel{\left( \begin{matrix} -x_1y_2 & 0 & 0 & x_2y_1 \\ x_2y_1 & -x_1y_2 & 0 & 0 \\ 0 & x_2y_1 & -x_1y_2 & 0 \\ 0 & 0 & x_2y_1 & -x_1y_2 \end{matrix} \right)}{\longrightarrow} & \bigoplus_{\chi\in \hat{G}} \mathcal{O}(\chi, -\chi) & \rightarrow & 0 \\ 
 & & \downarrow \\
 & & \Delta & \rightarrow & 0 
\end{array} 
\]

where the vertical map downwards takes the basis element $e_\chi$ for $\mathcal{O}(\chi,-\chi)$ to $z^{\Lambda(\overline{u})} \in \pi(M_{\Lambda(L)})$. 

\subsection{Weighted projective space $\PP(1,2)$}

For $\mathcal{L}, \mathcal{L}'$ line bundles on normal toric variety $X_\Sigma$, if $\mathcal{L}$ has sheaf of sections $\mathcal{O}(D)$ and $\mathcal{L}'$ has sheaf of sections $\mathcal{O}(E)$, then to compute

\[ Ext^i(\mathcal{O}_{X_\Sigma}(D), \mathcal{O}_{X_\Sigma}(E)) \]

we note that for $i=0$, 

\begin{align*}
    \mathcal{H}om(\mathcal{O}(D), \mathcal{O}(E)) \simeq \mathcal{O}(E-D)
\end{align*}

and $Ext^i(\mathcal{O}(D), \mathcal{O}(E)) = R\Gamma(\mathcal{O}(E-D)) = H^*(X_\Sigma, \mathcal{O}(E-D))$. \\

Now, for 

$$\PP(1,2) = \left[ \faktor{\C^2\setminus\{(0,0)\}}{\C^*} \right]$$

where $\C^*$ acts via $\lambda\cdot(z_1,z_2) = (\lambda z_1, \lambda z_2)$, the set \[ \mathcal{O}, \mathcal{O}(1), \mathcal{O}(2) \] is an exceptional collection. 

\begin{proof}
Note that 

\begin{align*}
H^*(\mathcal{O}(-1)) &= R\Gamma Hom(\mathcal{O}(2),\mathcal{O}(1))\\
&= R\Gamma Hom(\mathcal{O}(1),\mathcal{O}). 
\end{align*}

To find $H^*\mathcal{O}(n)$ on $\PP(1,2)$, we form a $\check{C}$ech cover on $\C^2$, noting that $z_1$ has degree $1$ and $z_2$ has degree $2$. On $U_1$, we localize at $z_1$ to get $\mathcal{O}(U_1)$. In Figure~\ref{fig:OU1.OU2} each lattice point is labeled by the total degree of monomial. On $U_2$, localizing at $z_2$ gives $\mathcal{O}(U_2)$. Here, lattice points in turqoise are again labeled by total degree of monomial. Level sets of total degree are highlighted in increments of 2. In the $\check{C}$ech cochain complex, $\check{C}^0(U_1 \cup U_2)$ is given by $\Gamma(U_1) \times \Gamma(U_2)$ in Figure~ \ref{fig:OU1.OU2}. Note that both $\Gamma(U_1)$ and $\Gamma(U_2)$ contain $(0,0)$. For \[H^*\left(\left[ \faktor{\C^2\setminus\{(0,0)\}}{\C^*}\right]\right)\] we consider elements only of degree $0$.\\[1cm]

\begin{center}
\begin{figure}[h]
\hspace{4cm} \includegraphics[trim = {2cm 5.5cm 3cm 2.5cm}, clip, width=.5\textwidth]{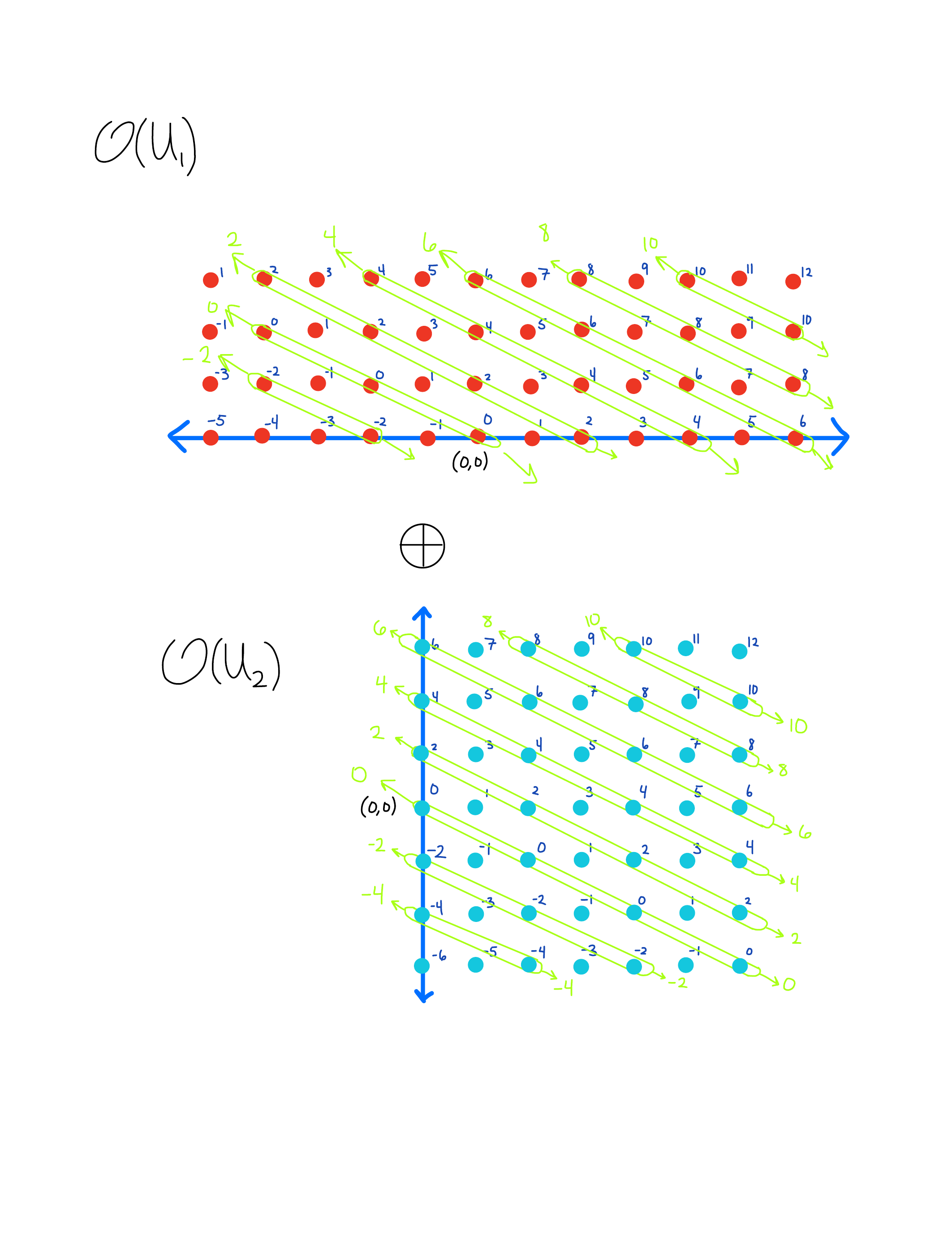}
\caption{$\mathcal{O}(U_1)$ and $\mathcal{O}(U_2)$ on $\PP(1,2)$}
\label{fig:OU1.OU2}
\end{figure} 
\end{center}


\begin{figure}[h]
\includegraphics[width=\textwidth]{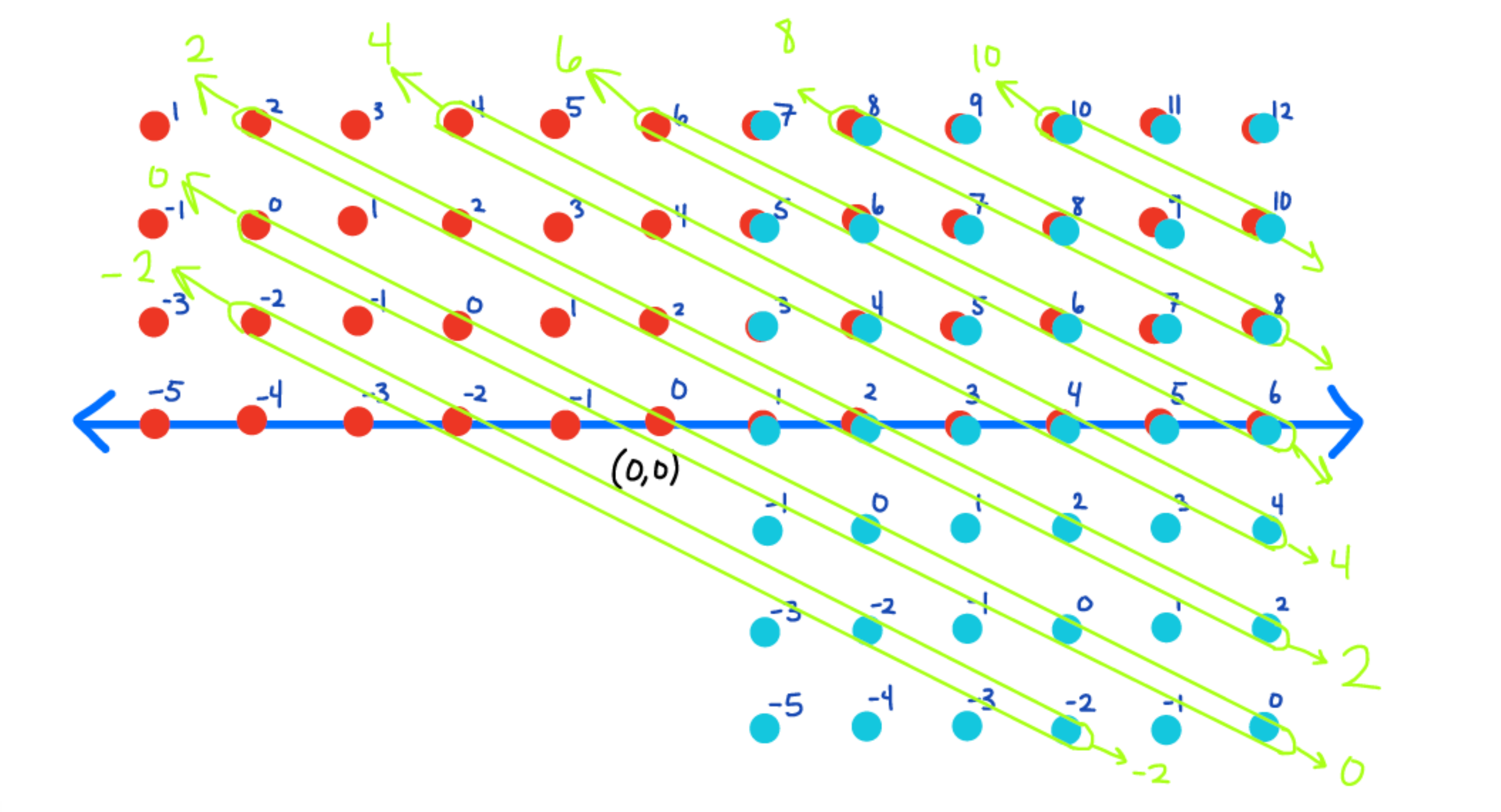}
\caption{$d^0(\check{C}^0)$ for $\mathcal{O}(-1)$ on $\PP(1,2)$}
\label{fig:C1}
\end{figure}

 To build the $\check{C}$ech complex for $\mathcal{O}(n)$ on $\PP(1,2)$, we map $\mathcal{O}(U_1)$ to $\mathcal{O}(U_2)$ using $d^0$, which shifts $\mathcal{O}(U_2)$ $n$ steps to the left. For $\mathcal{O}(-1)$, and $H^*(\mathcal{O}(-1))$, $d^0(\check{C}^0)$ and $\check{C}^1$ are given in Figure~\ref{fig:C1}. Here, $H^1$ is the kernel of $d^1$ by a dimension argument. Any lattice points labeled by both red and turquoise have elements in the kernel. But this does not occur for any points with degree 0, so $H^1(\PP(1,2),\mathcal{O}(-1))=0$. Here, $H^0$ is the cokernel: but all degree 0 elements lie in the image of $\check{C}^0$, since each degree 0 lattice point is labeled either red or turquoise. Therefore, \[ H^i(\PP(1,2), \mathcal{O}(-1)) = 0 \] for all $i$. \\[1cm]
 
 For $\mathcal{O}(-2)$, we would move the turquoise dots over to the right by $1$ in Figure~\ref{fig:C1} (a total of 2 steps to the right). However, we still have no lattice points of degree $0$ labeled by both red and turqoise lattice points, so $H^1(\PP(1,2),\mathcal{O}(-2))=0$. Similarly, all degree $0$ lattice points are labeled by either red or turqoise dots, so the cokernel is 0, and $H^0(\mathcal{O}(-2))=0$ on $\PP(1,2)$.

 \begin{figure}[h]
 \includegraphics[width=\textwidth]{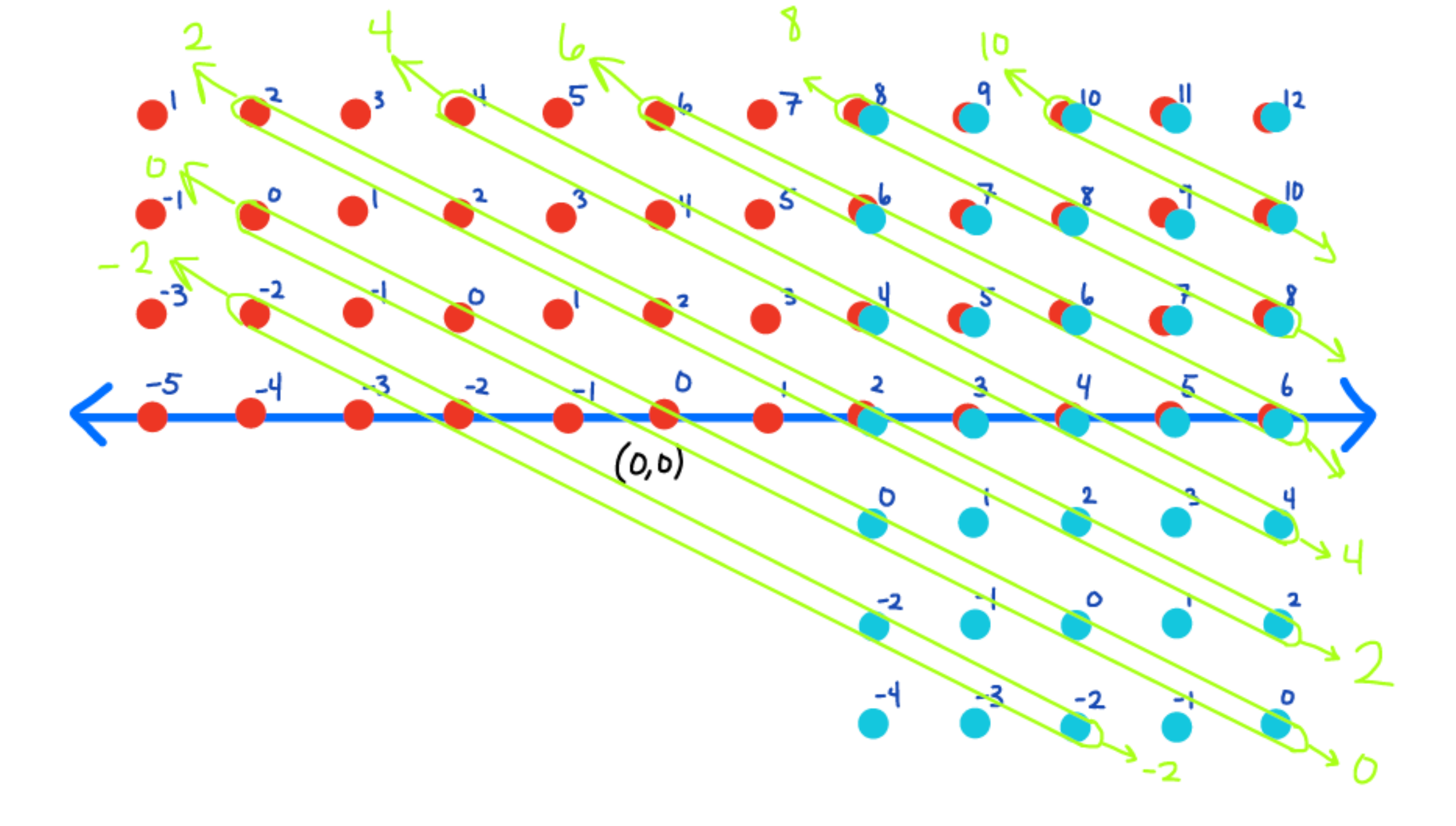}
 \caption{$d^0(\check{C}^0)$ for $\mathcal{O}(-2)$ on $\PP(1,2)$}
 \label{fig.}
 \end{figure} 
 
 Therefore

 \begin{align*}
 Ext(\mathcal{O}(2),\mathcal{O}(1)) \cong Ext(\mathcal{O}(1),\mathcal{O}) &\cong R\Gamma(\mathcal{O}(-1))\\
 &\cong H^*(\mathcal{O}(-1))  = 0 \end{align*}

 and

 \begin{align*}
Ext(\mathcal{O}(2),\mathcal{O}) = R\Gamma(\mathcal{O}(-2)) = H^*(\mathcal{O}(-2))=0. \end{align*} 

Note that $d^0(\check{C}^0)$ for $\mathcal{O}(-2)$ is given in Figure~\ref{fig.}. To see that $\mathcal{O},\mathcal{O}(1),\mathcal{O}(2)$ generate $D^b(\PP(a,b))$ for $a,b>0$ and coprime, we note exactness of \[ 0 \rightarrow \mathcal{O}(-a-b) \rightarrow \mathcal{O}(-a)\oplus\mathcal{O}(-b) \rightarrow \mathcal{O}\rightarrow 0 \] 

on $\PP(a,b)$. 
\end{proof}

\section{Resolving the diagonal for smooth toric variety $X_\Sigma$}

\subsection{Resolving the diagonal via deformation: $Bl_{[1:0:0]}(\PP^2)$ and deforming to $\mathcal{H}_L^\epsilon$ for projective toric varieties which are smooth and non-unimodular} 

\label{section3.8}

For toric varieties which are smooth and non-unimodular in the sense of Bayer-Popescu-Sturmfels, we'll need to deform the cellular complex $\mathcal{H}_L$ by an element $\epsilon\in \R^A$. We use $\epsilon$ to give transversal intersections at the vertices in $\mathcal{H}_L^\epsilon(0)$. Here, we deform the cellular complex $\mathcal{H}_L$ to $\mathcal{H}_L^\epsilon$ for $Bl_{[1:0:0]}\PP^2$ which is smooth, Fano, and unimodular as an illustrative example.\\[.2cm]

\begin{figure}[h]
\includegraphics[width=\textwidth]{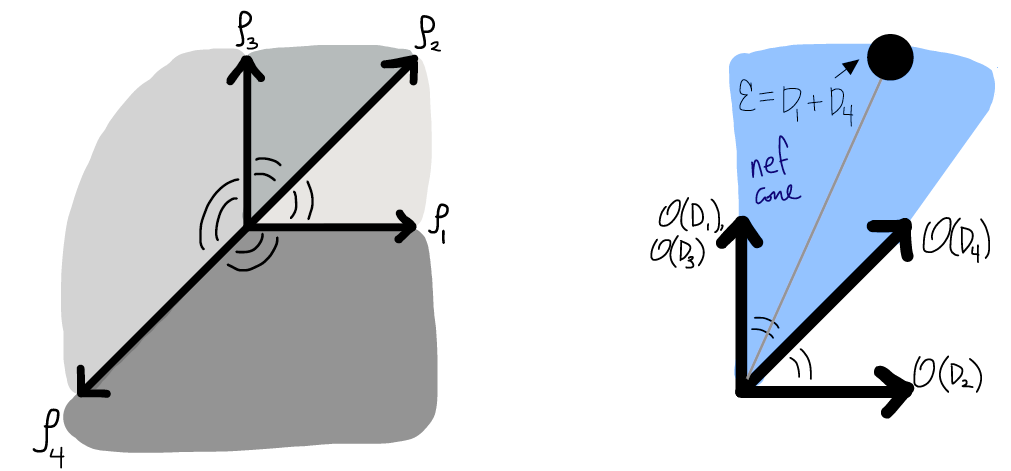}
\caption{The fan $\Sigma$ and the nef cone for $Bl_{[1:0:0]}\PP^2$ }
\label{fig:Bl0P2}
\end{figure}

Recall that a sublattice $L\subset \Z^n$ is called \textbf{unimodular} if $L$ is the image of an integer matrix $B$ with linearly independent columns, such that all maximal minors of $B$ lie in the set $\{0,1,-1\}$ \cite{bayer-popescu-sturmfels}. If we blow-up the torus-invariant point $[1:0:0]$ of $\PP^2$ (following Cox-Little-Schenck, I will call $X_\Sigma=Bl_{[1:0:0]}\PP^2$ the toric variety associated to the fan in $N_\R$ given in Figure~\ref{fig:Bl0P2}), then the lattice $L$ given by the image of $B$ in 

\[ 0 \rightarrow M \stackrel{B}{\longrightarrow} \Z^4 \stackrel{\pi}{\rightarrow}  Cl(X_\Sigma) \rightarrow 0 \] 

where $B$ given by $\left[\begin{matrix} 1 & 0 \\ 1 & 1 \\ 0 & 1 \\ -1 & -1 \end{matrix}\right]$ is unimodular.

Here, 

\[ Cl(X_\Sigma) \cong \faktor{\left<D_1, \dots, D_4\right>}{(D_1+D_2-D_4, D_2+D_3-D_4)}\]

Since $D_2$ corresponds to the exceptional divisor in the blow-up, $D_2$ has self-intersection number -1: \[ D_2\cdot D_2=-1\] while $D_1, D_3$, and $D_4$ have self-intersection number $1$, so the nef cone of $X_\Sigma$ is given by the span of $D_1$ and $D_4$. This description of the nef cone comes from a map

\[ (\R^2)^\vee \stackrel{\left[\begin{matrix} 0 & 1 & 0 & 1 \\ 1 & 0 & 1 & 1 \end{matrix}\right]}{\longleftarrow} \R^4 \]

corresponding to a fan in the first orthant of $\R^4$. Here, the element $D_1+D_4$ in the nef cone corresponds to $\left[\begin{matrix} \epsilon \\ 0 \\ 0 \\ \epsilon\end{matrix}\right] = \epsilon \left[\begin{matrix} 1 \\ 0 \\ 0 \\ 1 \end{matrix}\right] \in \R^4$. 

\subsubsection{Deforming the cellular complex $\mathcal{H}_L$}
Here, $X_\Sigma$ is assumed to be smooth.
If we consider $\R L = L\otimes_\Z \R \subset \R^4$ as before from the unimodular case, then the intersection of the infinite hyperplane arrangement $\mathcal{H}_L$ with $\R L$ will have more lattice points than just what appears in $L$ \cite{bayer-popescu-sturmfels}. That is, for $X_\Sigma$ smooth and non-unimodular,

\[ (\mathcal{H}_L \cap \R L) \setminus L \neq \emptyset. \] 

\begin{figure}[h]
\hspace*{3cm}  \includegraphics[trim = {4cm, 17cm, 4cm, 2cm}, clip, width=\textwidth]{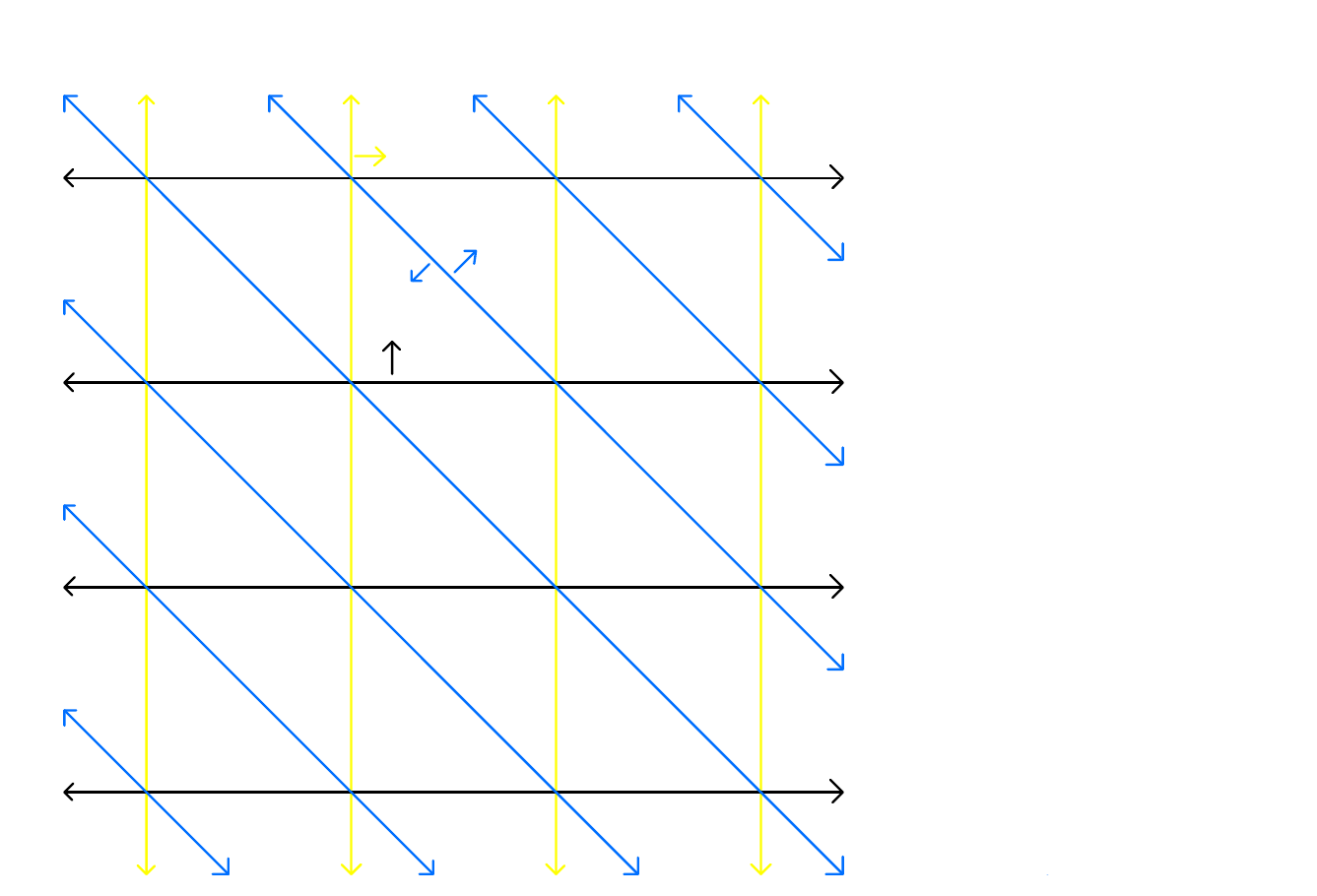}
\caption{$\R L$, showing $\R L \cap \mathcal{H}_L \subset \R^4$}
\label{fig:RLHL}
\end{figure}

$\R L\cap \mathcal{H}_L$ is given in Figure~\ref{fig:RLHL} for $Bl_{[1:0:0]}\PP^2$, which is unimodular (hence no extra vertices appear). Here, vertical lines in yellow correspond to integer values of $x_1$, and $x_1$ increases in value as we move to the right in the diagram. Horizontal lines in black correspond to integer values of $x_3$, which increase as we move up. Diagonal lines in blue give integer values in $x_2$ and $x_4$, which increase up and to the right, and down and to the left, respectively.

Since more than two hyperplanes intersect in $\R L$ at each point, the cellular complex $\mathcal{H}_L$ does not have transversal intersections at each vertex.

To remedy the lack of transversal intersection at each vertex $v \in \mathcal{H}_L(0)$, we adjust the infinite cellular complex $\mathcal{H}_L$ by deforming the hyperplane arrangement in a manner governed by our choice of $\epsilon$ in the nef cone, given by

\[ x \left[ \begin{matrix} 1 \\ 1 \\ 0 \\ -1 \end{matrix}\right] + y\left[\begin{matrix} 0 \\ 1 \\ 1 \\ -1 \end{matrix}\right] + \left[\begin{matrix} \epsilon \\ 0 \\ 0 \\ \epsilon \end{matrix}\right] = 0 \] 

corresponding to row $1$ and row $2$ of $B$, giving a basis for $L\subset \R^4$. 

\begin{figure}[h]
\includegraphics[trim = {4cm, 15cm, 6cm, 3cm}, clip, width=\textwidth]{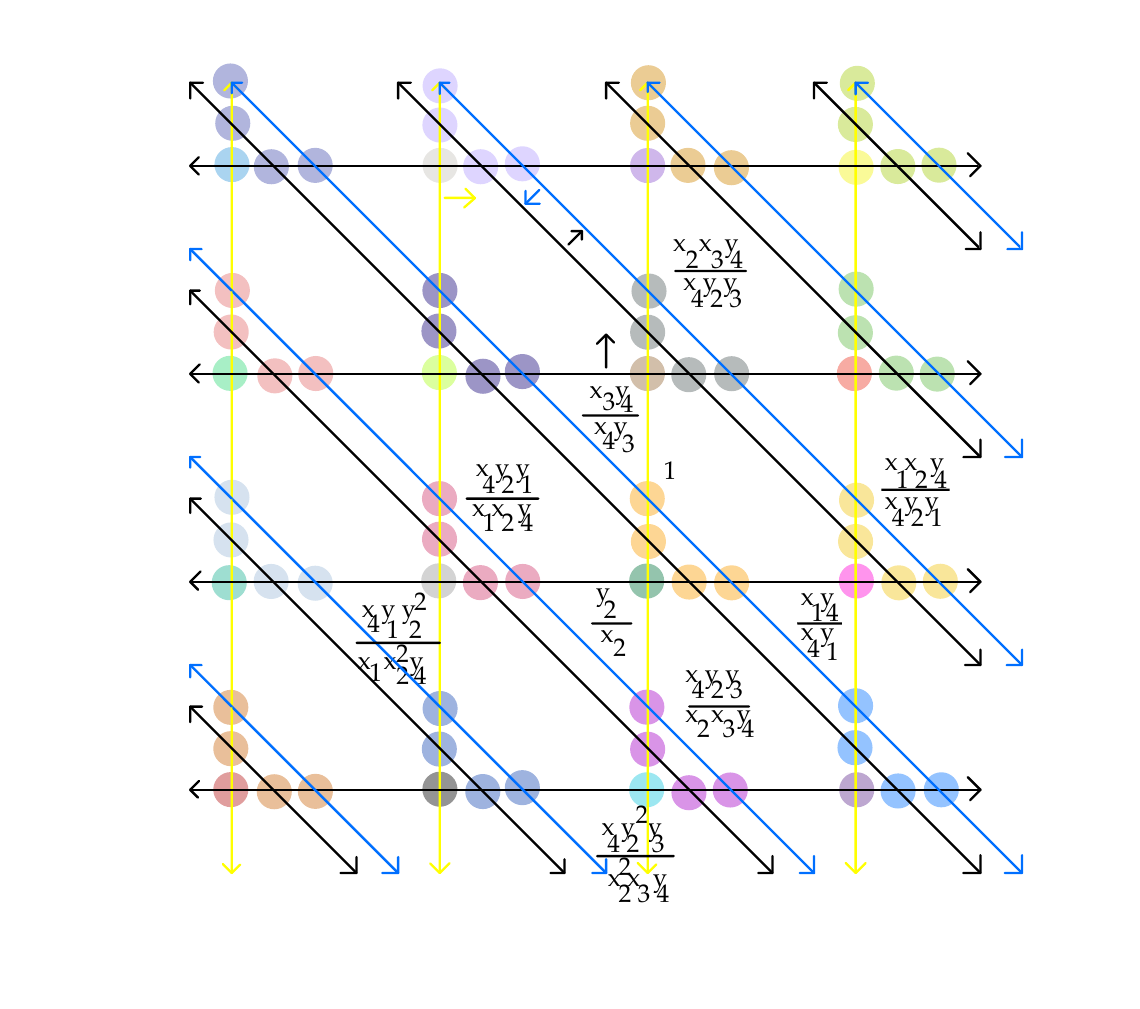}
\caption{$\mathcal{H}_L^\epsilon$ for $Bl_{[1:0:0]}\PP^2$ deformed by $\mathcal{O}(D_1)+\mathcal{O}(D_4)$}
\label{fig:HLepsilon}
\end{figure}

This gives the deformed hyperplane arrangement $\mathcal{H}_L^\epsilon$ corresponding to the shifts $e_1 + \epsilon \in \Z, e_4 + \epsilon \in \Z$ for $0<\epsilon \ll 1$. We construct the Hull resolution $(\mathcal{F}^\bullet_{\mathcal{H}_L^\epsilon}, \partial^\epsilon)$ below. The deformed complex $\mathcal{H}_L^\epsilon$ is given in Figure~\ref{fig:HLepsilon}, together with the monomial labelings given by taking the vertex labeling associated to the floor function for each coordinate of a given point in $\mathcal{H}_L^\epsilon \subset \R^4$. Vertices with common color labeling are given the same monomial vertex labeling in Figure~\ref{fig:HLepsilon}. The arrows off of hyperplanes indicate direction of increase of $x_i$ for $1\leq i \leq 4$, as in Figure~\ref{fig:RLHL}. 

Here, we emphasize that $\R L$ is a plane in $\R^4$ along which we are considering the intersections with the hyperplane arrangement $\mathcal{H}_L = \{ x_i \in \Z \text{ }|\text{ } 1\leq i \leq 4 \}$. 

\begin{figure}[h]
\hspace{3cm} \includegraphics[trim = {6cm 14cm 3cm 2cm}, clip,width=.6\textwidth]{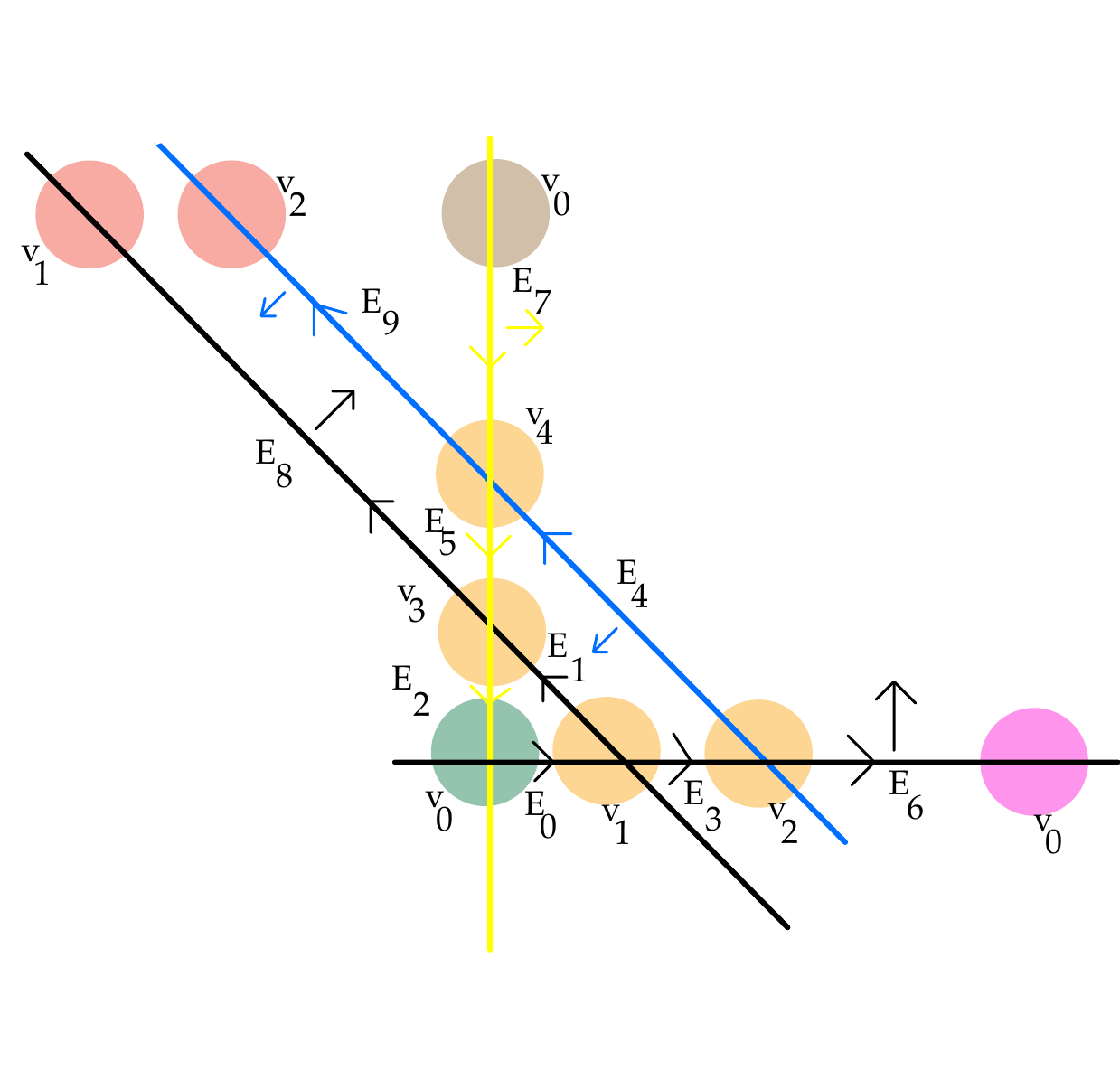}
\caption{The quotient cellular complex $\faktor{\mathcal{H}_L^\epsilon}{L}$ with monomial labelings from $M_{\Lambda(L)}$}
\label{fig:HLeQuotL}
\end{figure}

Taking the quotient of $\mathcal{H}_L^\epsilon$ by $L$ and giving monomial labels from $M_{\Lambda(L)}$ given by the integer floor function in each component are given in Figure~\ref{fig:HLeQuotL}.\\[2cm]
\newpage

Here, the differential 

\[ \partial^\epsilon_1: \bigoplus_{e \in \mathcal{H}_L^\epsilon(1)}S \rightarrow \bigoplus_{v\in \mathcal{H}_L^\epsilon(0)} S \] is given by 

\[\begin{array}{c|cccccccccc|}
v_4 & 0 & 0 & 0 & 0 & 1 & -1 & 0 & x_3y_4 & 0 & -x_3y_1\\ v_3 & 0 & 1 & -y_2 & 0 & 0 & 1 & 0 & 0 & -x_3y_1 & 0 \\ 
v_2 & 0 & 0 & 0 & 1 & -1 & 0 & -x_1y_4 & 0 & 0 & x_1y_3\\
v_1 & y_2 & -1 & 0 & -1 & 0 & 0 & 0 & 0 & x_1y_3 & 0 \\
v_0 & -x_2 & 0 & x_2 & 0 & 0 & 0 & x_4y_1 & -x_4y_3 & 0 & 0 \\
\hline
 & E_0 & E_1 & E_2 & E_3 & E_4 & E_5 & E_6 & E_7 & E_8 & E_9 \end{array} \]

 subject to the labeling in Figure~\ref{fig:HLeQuotL}. The color labeling of vertices and edges in Figure~\ref{fig:HLeQuotL} corresponds to the color labeling in Figure~\ref{fig:HLepsilon}. Exactness away from homological index $i=0$ follows as before for the unimodular case by the same convexity argument used in Theorem 3.1 in Bayer-Popescu-Sturmfels\cite{bayer-popescu-sturmfels}. We investigate the cokernel of $\partial^\epsilon_0$ in the following section.

\subsection{Cokernel of $\partial^\epsilon_0$}
\label{nefcone}

Here, we consider the map

\[ (\mathcal{F}_{\mathcal{H}_L^\epsilon}, \partial^\epsilon) \stackrel{f}{\rightarrow} S_{\prod x_i}\]

given by $\bigoplus_{v\in \mathcal{H}_L^\epsilon}S(-m_F) \rightarrow S_{\prod x_i}$ where a vertex in $\mathcal{H}_L^\epsilon$ maps to its vertex monomial label in $M_{\Lambda(L)}$.

\begin{figure}[h]
\includegraphics[width=\textwidth]{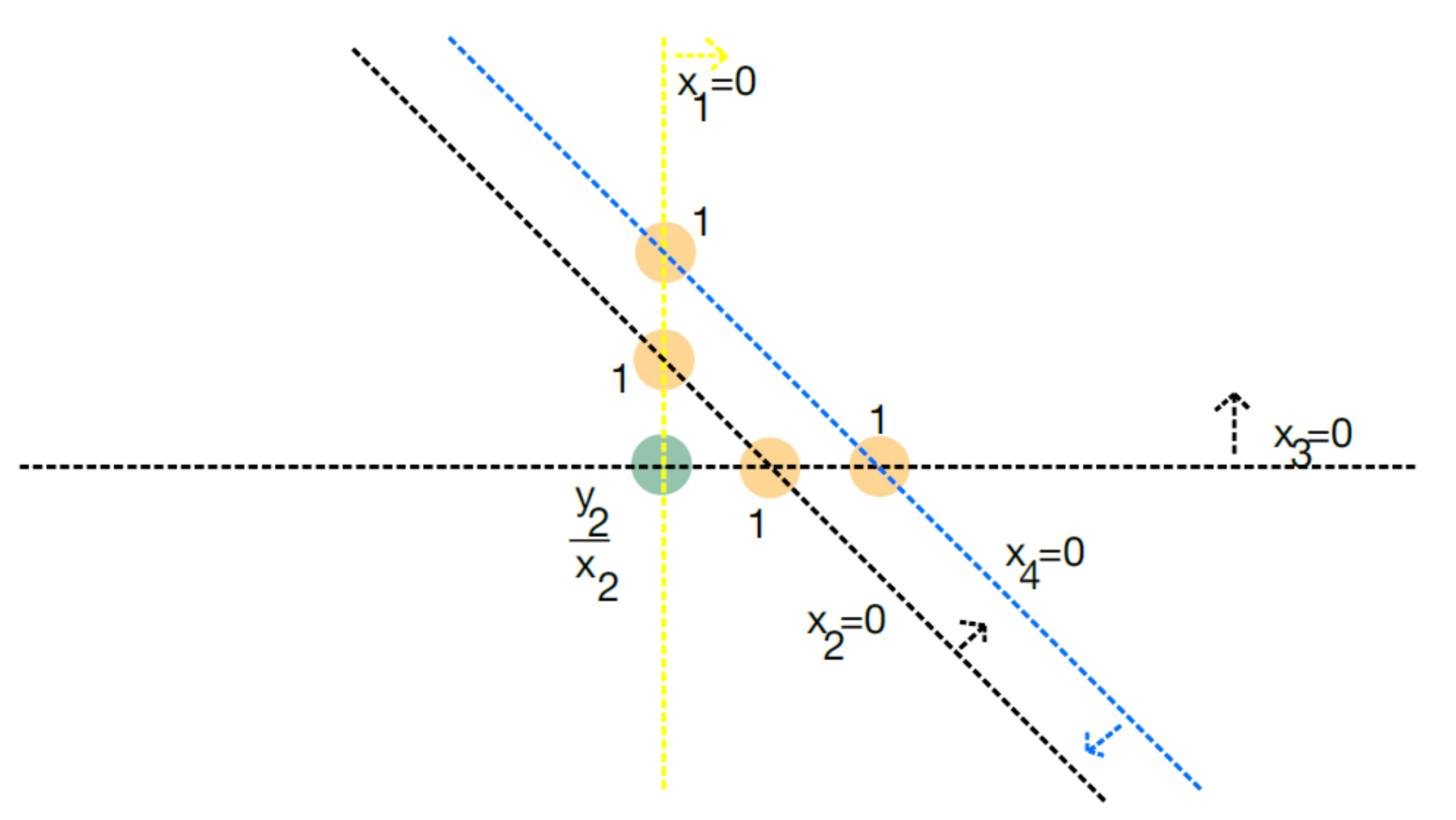}
\caption{Monomial labels on vertices $v\in \mathcal{H}_L^\epsilon(0)$ modulo $L$}
\label{fig: HLL.12.7}
\end{figure}

For $Bl_{[1:0:0]}\PP^2$, we have that $L = \left< \left[\begin{matrix} 1 \\ 1 \\ 0 \\ -1\end{matrix}\right] , \left[\begin{matrix} 0 \\ 1 \\ 1 \\ -1 \end{matrix}\right]\right>$, \\

\begin{align*}
I_L &= \left< x^{v_+}-x^{v_-} \text{ }|\text{ } v\in L \right> \\
 &= \left< x_1x_2-x_4, x_2x_3-x_4, x_1-x_3\right>, \end{align*}

 \begin{align*}
J_L &= \left< x^uy^v - x^vy^u \text{ }|\text{ }u-v\in L\right>\\
&= \left<x_1x_2y_4 - x_4y_1y_2, x_2x_3y_4 - x_4y_2y_3, x_1y_3-x_3y_1\right> \end{align*}

and $M_{\Lambda(L)}$ as an $S-$submodule of $T$ has the infinite generating set $\{ x^uy^{-u}\text{ }|\text{ }u\in L \}$. In the quotient by $L$, the cellular complex $\mathcal{H}_L^\epsilon$ carries monomial labels given in Figure~\ref{fig: HLL.12.7}. When we consider monomial labels on vertices in the quotient $\mathcal{H}_L^\epsilon/L$, we only need to consider the monomial labels on vertices given in Figure~\ref{fig: HLL.12.7} due to which vertices in $\mathcal{H}_L^\epsilon(0)$ carry the same monomial label. \\

Here, $L\cap \N^n = \{0\}$ implies that $M_{\Lambda(L)}\ni 1$, so that $M_{\Lambda(L)} \subset Im(f)$ above, by considering the action of $L$ on the quotient from Figure~\ref{fig: HLL.12.7}. Modulo the action of $L$, there is one additional monomial $\frac{y_2}{x_2}$ in $Im(f)$ which is not contained in $M_{\Lambda(L)}$. Note that since 

\[ \{ \frac{x_1x_2 y_4}{x_4y_1y_2}, \frac{x_2x_3y_4}{x_4y_2y_3}, \frac{x_1y_3}{x_3y_1} \} \subset M_{\Lambda(L)},\]

we also have

\[ \{\frac{x_4y_1y_2}{x_1x_2y_4}, \frac{x_4y_2y_3}{x_2x_3y_4}, \frac{x_3y_1}{x_1y_3} \} \subset M_{\Lambda(L)}. \] 

Now for $X_\Sigma = Bl_{[1:0:0]}\PP^2$,

\begin{align*}
I_{irr} &= \left< x_3x_4, x_1x_4, x_1x_2, x_2x_3\right> \\
&= \left<x_1, x_3\right> \cap \left<x_2,x_4\right> \end{align*}

and $I_{irr}$ for $X_\Sigma \times X_\Sigma$ is given by \[ \left<x_1, x_3\right>\cap\left<x_2,x_4\right> \cap \left<y_1,y_3\right> \cap \left<y_2, y_4\right>. \] 

Now, for any monomial $q$ in $I_{irr}$ for $X_\Sigma \times X_\Sigma$, we must have that either $x_2$ or $x_4$ divides $q$. If $x_2$ divides $q$, then

\[ x_2 \left( \frac{y_2}{x_2} \right) = y_2 \in M_{\Lambda(L)} \] 

since $M_{\Lambda(L)} \ni 1$ shows that $q \cdot \frac{y_2}{x_2} \in M_{\Lambda(L)}$. If $x_4 | q$, then 

\begin{align*}
y_3x_4 \cdot \left( \frac{y_2}{x_2}\right) &= \left(\frac{x_4y_2y_3}{x_2x_3y_4}\right) \cdot x_3y_4 \in M_{\Lambda(L)} \end{align*}

shows that $q \cdot \left(\frac{y_2}{x_2}\right)$ is in the image of the action of $S$ on $q\cdot \left( \frac{y_2}{x_2}\right)$. Therefore, the cokernel of the inclusion of $M_{\Lambda(L)} \hookrightarrow Im(f)$ is torsion with respect to $I_{irr}$ for $X_\Sigma \times X_\Sigma$. In general, we say that an $R-$module $M$ is $I$ torsion iff there exists $k\in \Z_+$ such that $I^k M = 0$. Here, $k=1$ suffices to show that the cokernel of the inclusion of $M_{\Lambda(L)} \hookrightarrow Im(f)$ is torsion with respect to $I_{irr}$ for $X_\Sigma \times X_\Sigma$. \\[2cm]

\subsection{Moving around in the effective cone}

\begin{figure}[h]
\includegraphics[width=\textwidth]{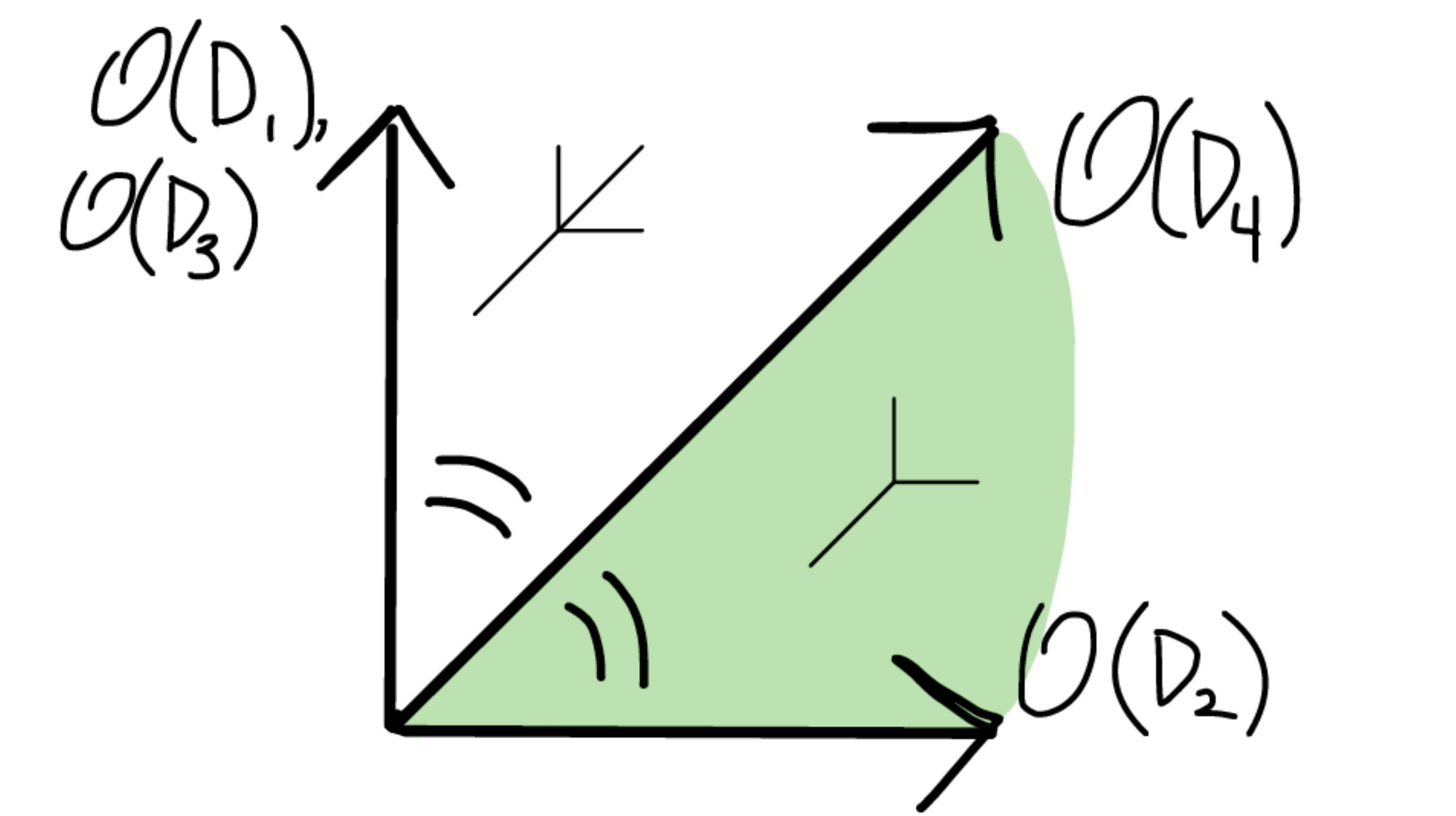}
\caption{Chambers of the Effective Cone}
\label{fig:effectivecone}
\end{figure}

Here we consider $\text{coker }(M_{\Lambda(L)} \hookrightarrow Im(f))$
for $\epsilon$ in the interior of the other chamber of the effective cone, spanned by $\mathcal{O}(D_2)$ and $\mathcal{O}(D_4)$. (Again, $D_2$ corresponds to the exceptional divisor of the blow-up of $Bl_{[1:0:0]}\PP^2$.) Our previous discussion of $\text{coker}(M_{\Lambda(L)}\hookrightarrow Im(f))$ from Section~\ref{nefcone}  took place inside the nef cone of $N^1(X_\Sigma)$, with $N^1(X_\Sigma)\cong Cl(X_\Sigma)\otimes_\Z \R$ since $X_\Sigma$ is smooth. More generally, we can ask about the cokernel of the inclusion of $M_{\Lambda(L)}\hookrightarrow Im(f)$ inside of the other chamber of the effective cone, given in Figure~\ref{fig:effectivecone}.

Inside of this cone spanned by $\mathcal{O}(D_2)$ and $\mathcal{O}(D_4)$ in $Cl(X_\Sigma)\otimes_\Z \R$, we have a corresponding fan for which the ray $\rho_2$ is removed. This gives the fan for $\PP^2$, though we retain the information that $\rho_2$ no longer lives in any maximal cone so that the irrelevant ideal for $X_\Sigma$ becomes

\[ I_{irr} = \left<x_1,x_3,x_4\right> \cap (x_2) \]

and the irrelevant ideal for $X_\Sigma \times X_\Sigma$ in this case is

\[ I_{irr} = \left<x_1,x_3,x_4\right> \cap (x_2)\cap \left<y_1, y_3, y_4\right> \cap (y_2). \] 

\begin{figure}
\includegraphics[width=.8\textwidth]{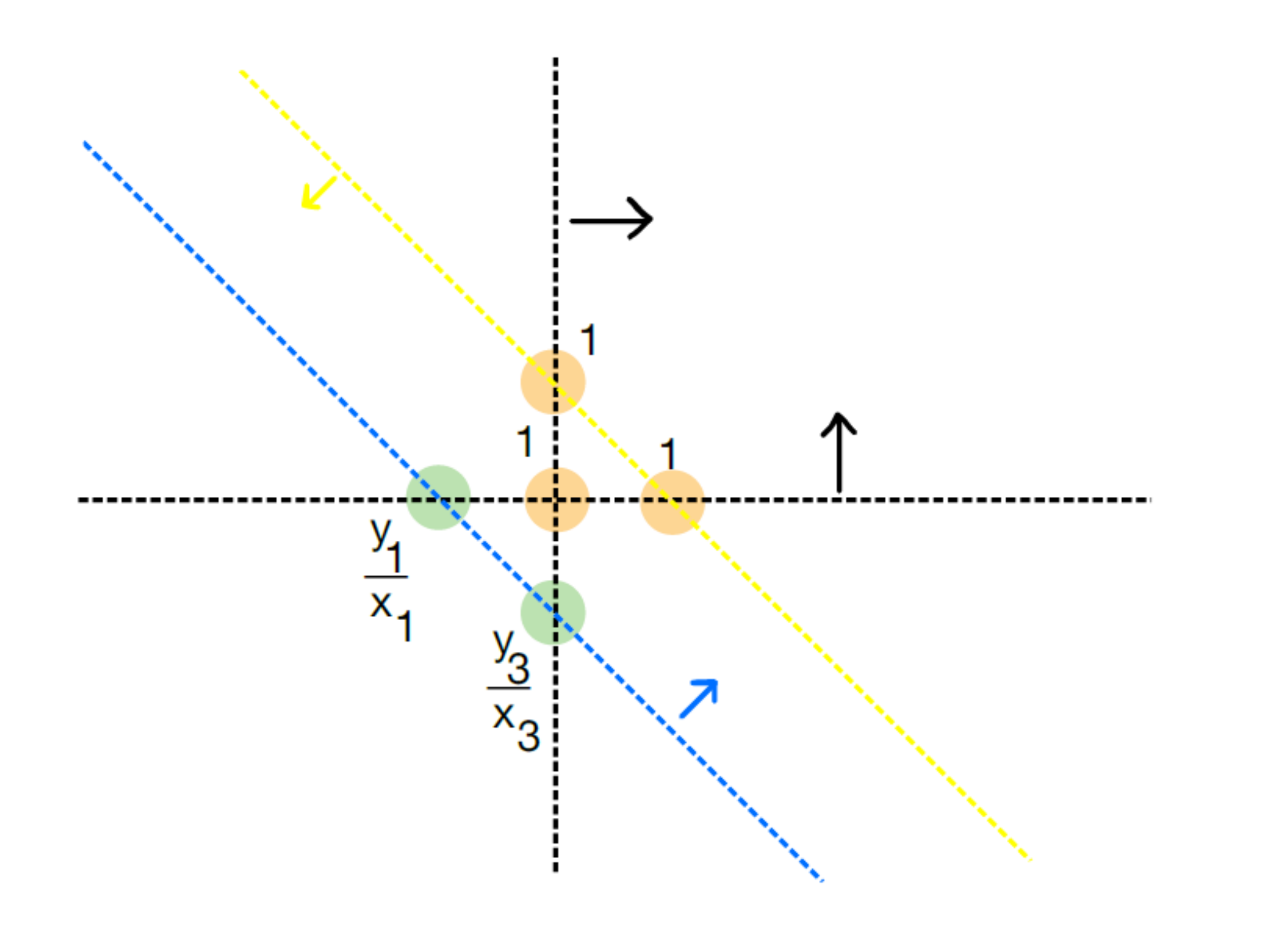}
\caption{Deformed complex $\mathcal{H}_L^\epsilon$ for $\epsilon$ corresponding to $\mathcal{O}(D_2) + \mathcal{O}(D_4)$}
\label{fig: def2}
\end{figure}

Taking $\epsilon = \left[ \begin{matrix} 0 \\ \epsilon \\ 0 \\ \epsilon \end{matrix}\right]$ corresponding to $\mathcal{O}(D_2)+\mathcal{O}(D_4)$ in the interior of the chamber of the effective cone spanned by $\mathcal{O}(D_2)$ and $\mathcal{O}(D_4)$, we have the deformed cellular complex $\mathcal{H_L^\epsilon}$ whose quotient by $L$ is given in Figure~\ref{fig: def2}. Here, $x_2=0$ and $x_4=0$ are labeled in blue and yellow, respectively, since those planes are deformed in $\R L$. 
In this case, $I_{irr}$ for $X_\Sigma \times X_\Sigma$ is given above and $\text{coker}(M_{\Lambda(L)}) \ni \frac{y_1}{x_1}, \frac{y_3}{x_3}$ and is spanned by translates of these elements by $L$. Note that in this case, at least one of $x_1, x_3$, or $x_4$ must divide any monomial $q$ in the irrelevant ideal $I_{irr}$ of $X_\Sigma\times X_\Sigma$. \\[.2cm]

For $L$-translates of $\frac{y_1}{x_1}$, we note that if $x_1$ divides $q$, then 

\begin{align*}
x_1 \cdot (\frac{y_1}{x_1})=y_1 \in M_{\Lambda(L)}\end{align*} since $1\in M_{\Lambda(L)}$, and $M_{\Lambda(L)}$ is an $S-$submodule of $T= k[x_1^{\frac{+}{}}, \dots, x_n^{\frac{+}{}}, y_1^{\frac{+}{}}, \dots, y_n^{\frac{+}{}}]$. If $x_3$ divides $q$, then 

\begin{align*}
x_3 \cdot \frac{y_1}{x_1} &= y_3 \cdot \left(\frac{x_3y_1}{x_1y_3}\right) \in M_{\Lambda(L)}. \end{align*}

If $x_4$ divides $q$, then

\begin{align*}
x_4y_2 \left(\frac{y_1}{x_1}\right) &= x_2y_4\left(\frac{x_4y_1y_2}{x_1x_2y_4}\right) \in M_{\Lambda(L)}\end{align*}

so that $\overline{\frac{y_1}{x_1}} = \overline{0} \in \text{coker}(M_{\Lambda(L)}\hookrightarrow Im(f) )$.

Next, for ($L-$translates of) $\frac{y_3}{x_3}$, we again have that at least one of $x_1, x_3$, or $x_4$ divides any monomial $q$ in $I_{irr}$, the irrelevant ideal for $X_\Sigma\times X_\Sigma$. If $x_1$ divides $q$, then 

\begin{align*}
\frac{y_3}{x_3}\cdot x_1 &= \left(\frac{x_1y_3}{x_3y_1}\right)\cdot y_1 \in M_{\Lambda(L)}\end{align*}

If $x_3$ divides $q$, then 

\begin{align*} \frac{y_3}{x_3}\cdot x_3 &= y_3\in M_{\Lambda(L)} \end{align*} since $1 \in M_{\Lambda(L)}$. Lastly, if $x_4 | q$, then \begin{align*} \left(\frac{y_3}{x_3}\right)x_4y_2 &= \left(\frac{x_4y_2y_3}{x_2x_3y_4}\right)x_2y_4 \in M_{\Lambda(L)}. \end{align*}

This shows that $\overline{\frac{y_3}{x_3}} = 0 \in \text{coker}(M_{\Lambda(L)}\hookrightarrow Im(f))$. Hence,

\[ \text{coker}(M_{\Lambda(L)}\hookrightarrow Im(f)) = 0 \] 

modulo $I_{irr}$ for $X_\Sigma \times X_\Sigma$ for $\epsilon = \mathcal{O}(D_2)+\mathcal{O}(D_4)$ in the interior of the chamber spanned by $\mathcal{O}(D_2)$ and $\mathcal{O}(D_4)$ in the effective cone of the secondary fan. 

\label{theorem3.8.3}
\subsection{General argument for $\text{coker}(M_{\Lambda(L)}\hookrightarrow Im(f))$ for $\epsilon$ in the effective cone}

Given $X_\Sigma$ smooth and projective with no torus factors, we have the following exact sequence:

\[
\begin{array}{ccccccccc}
0 &\leftarrow &Cl(T)_\R &\leftarrow &\R^A & \leftarrow & L_\R &\leftarrow &0 \\
 & &  & & \rotatebox[origin=c]{90}{$\in$}\\
 & & & & \epsilon\end{array}\]

 where $\epsilon\mapsto \overline{\epsilon}$ in the effective cone of $\mathcal{F}_{\Sigma(A)}$. Now $\overline{\epsilon}$ gives a convex function $g_\epsilon: N_\R \rightarrow \R$ which is linear on all cones of a corresponding fan $\Sigma_\epsilon$ such that the maximal domains of linearity correspond to maximal cones in $\Sigma_\epsilon$. So $g_\epsilon$ gives the fan $\Sigma_\epsilon$ with corresponding irrelevant ideal $I_{irr}$ for $\Sigma_\epsilon\times\Sigma_\epsilon$ given by \[ \left< \prod_{i\in I_\emptyset} x_i \prod_{\rho\not\in\sigma} x_\rho \text{ }|\text{ } \sigma \in \Sigma_{max} \right> \cap \left< \prod_{i\in I_\emptyset} y_i \prod_{\rho\not\in \sigma} y_\rho \text{ } \sigma \in \Sigma_{max}\right> \] 

 Continuing with the map $(\mathcal{F}_{\mathcal{H}_L^\epsilon}^\bullet, \partial_\epsilon^\bullet) \stackrel{f}{\rightarrow} S_{\prod x_i}$ and our consideration of $coker(M_{\Lambda(L)} \hookrightarrow Im(f))$ modulo $I_{irr}$ for $\Sigma_\epsilon \times \Sigma_{\epsilon}$, we consider $p\in \R^{\mathcal{A}}$ such that $p \mapsto \overline{p}\in\mathcal{H}_L^\epsilon(0)$ near $\overline{0}$ in $\R L$. We have previously referred to $\overline{p}\in \R L$ as $p$. Here, we take 

 \[ \overline{p} = \bigcap_{\rho\in \mathcal{B}_p} \{x_\rho = 0 \} \]
 as a definition of $\mathcal{B}_p$. Now transversality of $\mathcal{H}_L^\epsilon$ implies that $\mathcal{B}_p$ gives a basis for $\R L$. Since $p$ is near $\vec{0}$ in $\R^\mathcal{A}$, the monomial label $m_p$ is squarefree, and $\epsilon \in \R^\mathcal{A}_{\geq 0}$ implies that $m_p$ can be written $\frac{y_{i_1} \cdots y_{i_\ell}}{x_{i_1}\cdots x_{i_\ell}}$. 
 
 By the definition of $\mathcal{H}_L^\epsilon$, the monomial label on $p$ is given by \[ m_p = \frac{x^{\vec{u}}}{y^{\vec{u}}} \text{ }|\text{ } \vec{u} = \lfloor p \rfloor. \] We now adjust $g_\epsilon$ by subtracting a linear functional $\ell_p$ determined by $\ell_p(u_\rho)=\epsilon_\rho$ for $\rho\in \mathcal{B}_p$. $\ell_p$ is determined by this information, since $\mathcal{B}_p$ gives a basis for $\R L$ and for $N$. That is, we construct \[ \tilde{g}_\epsilon := g_\epsilon - \ell_p \] so that we have \[ p = \tilde{g}_\epsilon(\{u_\rho\}) \] 

except for coordinates in $I_\emptyset$. Since any variable corresponding to an element of $I_\emptyset$ always appears in the support of any monomial in the irrelevant ideal, we may assume without loss of generality that $p = \tilde{g}_\epsilon$. To show that $coker (M_{\Lambda(L)}\hookrightarrow Im(f)) =0$ modulo $I_{irr}$, let $\sigma$ be any maximal cone in $\Sigma_{max}$. From convexity of $g_\epsilon$, it follows that $\tilde{g}_\epsilon$ is also convex. This gives a subdivision of maximal cones in $\Sigma$ such that either:

\label{pnear0}
\begin{itemize}
\item[(i)] $\sigma \subseteq \tilde{g}_\epsilon^{-1}(\R_{\geq 0})$\\
\item[(ii)] $\sigma \subseteq \tilde{g}_\epsilon^{-1}(\R_<0)$\\
\item[(iii)] Neither (i) nor (ii) holds \end{itemize}

In case (i), since we only need to worry about what's not in $I_\emptyset$, we can suppose that $g_\epsilon(u_\rho)=\epsilon_\rho$ so that $p = \tilde{g}_\epsilon(\{u_\rho\})$. If we choose $\epsilon$ to be small so that $|\epsilon|<<1$, then $\ell_p$ remains small as a continuous function of $\epsilon$ which gives that the coordinates of $p\in (-1,1)$ and the floor function of all coordinates is in the set $\{-1,0\}$. Therefore, the set of indices $\{i_j\}$ for which $x_{i_j}$ is in the support of $m_p$ are exactly the coordinates for which the floor function returns $-1$. Since $x^{\hat{\sigma}} = \prod_{\rho\not\in\sigma} x_\rho$ is a generator of $I_{irr}$ (ignoring $I_\emptyset$), we have that all variables cooresponding to rays of $\sigma$ do not lie in the support of $m_p$, so that $x^{\hat{\sigma}}$ can clear denominators to give $x^{\hat{\sigma}}\cdot m_p \in M_{\Lambda(L)}$.\\

In case (ii), we have that $\sigma \subseteq \tilde{g}_\epsilon^{-1}(\R_{<0})$, $m_p = \frac{y_{i_1}\cdots y_{i_\ell}}{x_{i_1} \cdots x_{i_\ell}}$ and $x^{\hat{\sigma}} = \prod_{\rho\not\in \sigma} x_\rho$. For $p$ near $\vec{0}$ we show that there exists $k>0$ such that $I_{irr}^k \cdot m_p \in M_{\Lambda(L)}$ and $coker(M_{\Lambda(L)}\hookrightarrow Im(f)) = 0$ by adding in the assumption that $X_\Sigma$ is Fano. If $X_\Sigma$ is Fano, then the anticanonical divisor $D = \sum_{\rho\in\Sigma(1)} D_\rho $ is very ample so that we can define $v\in L$ to be $v(u_\rho)=-1$ for $\rho\in \sigma$ and away from $\sigma$, $v(\rho') > -1$ for $\rho'\not\in \sigma$. Since $v\in L$ is an element of the lattice $L$, $v(u_\rho) \in \Z$ for all $\rho\in \Sigma(1)$ implies $v(\rho')\geq 0$ for all $\rho'\not\in \sigma$. 

 Now $\frac{x^{\vec{v}}}{y^{\vec{v}}} \in M_{\Lambda(L)}$ and we need only consider variables in the support of $m_p$ for which $\rho\in \sigma$ from our previous description of $I_{irr}$. Here, there exists some $\vec{w}\in \N^\mathcal{A}$ with $y^{\vec{w}}\in S$ and some $\{s_\rho\} \subset \N$ so that 

\begin{align*}
\prod_{\rho\not\in \sigma} x_\rho^{s_\rho} \cdot m_p &= y^w \cdot \frac{x^{\vec{v}}}{y^{\vec{v}}} \in M_{\Lambda(L)}
\end{align*}

and \[ (x^{\hat{\sigma}})^{\text{max}\{s_\rho\}} \cdot m_p = q \cdot \prod_{\rho\not\in \sigma} x_\rho^{s_\rho} m_p \in M_{\Lambda(L)} \] for some $q\in S$. Here we use the fact that $|\mathcal{A}|<\infty$ to take a maximum, and the fact that $I_{irr}$ has a finite generating set since $S$ is Noetherian.  \\

We generalize statements (i)-(iii) in the next subsection without the Fano assumption and therefore omit consideration of (iii) here. 

\subsubsection{For $p\in \mathcal{H}_L^\epsilon(0)$ away from $\vec{0}$}

For $p$ away from $\vec{0}$ and $X_\Sigma$ smooth (not assumed to be Fano here), we make use of the following lemmas:\\

\label{lemmata}

\textbf{Lemma 1}: Given $p\in \mathcal{H}_L^\epsilon(0)$ and $\sigma_1, \sigma_2 \in \Sigma(n)$, there exists $k>0$ with $(x^{\hat{\sigma}_1}y^{\hat{\sigma}_2})^k \cdot m_p \in M_{\Lambda(L)}$ iff there exists $v\in L$ and an $\ell>0$ such that $(x^{\hat{\sigma}_1}y^{\hat{\sigma}_2})^\ell \cdot m_{p+v} \in M_{\Lambda(L)}$. \\

\begin{proof} Any $v\in L$ has integer coordinates so that $\lfloor v \rfloor = v$ and $m_v = \frac{x^{\lfloor v \rfloor}}{y^{\lfloor v \rfloor }} = \frac{ x^v}{y^v}$. This implies that adding $v\in L$ to $p$ shifts all coordinates by only integer amounts, so that no cancellation occurs when we take a floor function: \[ \lfloor p+v \rfloor = \lfloor p \rfloor + v \]  

In particular, we note that if both $p$ and $v$ are in $L$, then $m_{p+v}=m_p \cdot m_v$ since $p, v \in L \implies p+v\in L$. So if there exists $k>0$ with $(x^{\hat{\sigma}_1}y^{\hat{\sigma}_2})^k \cdot m_p \in M_{\Lambda(L)}$, then $(x^{\hat{\sigma}_1}y^{\hat{\sigma}_2})^k\cdot m_{p+v} = (x^{\hat{\sigma}_1}y^{\hat{\sigma}_2})^k\cdot m_p m_v \in M_{\Lambda(L)}$ since $M_{\Lambda(L)}$ is closed under the action of $S[L]$. Now closure of $M_{\Lambda(L)}$ under the action of $S[L]$ also proves the reverse direction. \end{proof}


For the following \textbf{Claim}, we require the \textbf{Separation Lemma}\cite{fulton}, which separates convex sets by a hyperplane:

\textbf{Separation Lemma}: If $\sigma$ and $\sigma$' are convex polyhedral cones whose intersection $\tau$ is a face of each, then there is a $u$ in $\sigma^\vee \cap (-\sigma')^\vee$ with \[ \tau = \sigma \cap u^\perp = \sigma' \cap u^\perp. \]

Generalizing the previous cases (i)-(iii) from Section~\ref{pnear0}, we make the following claim.

\textbf{Claim}: Given $p\in \mathcal{H}_L^\epsilon(0)$, for all $\sigma_1, \sigma_2\in \Sigma(n),$ there exists $u\in L, k\in \mathbb{N}$ and $\alpha,\beta\in \mathbb{N}^{\Sigma(1)}$ such that $(x^{\hat{\sigma}_1}y^{\hat{\sigma}_2})^km_p = x^\alpha y^\beta m_u \in M_{\Lambda(L)}$. \\

\begin{proof}
The \textbf{Claim} holds iff there exists $u\in L, k\in \mathbb{N}$ and $\alpha,\beta\in \mathbb{N}^{\Sigma(1)}$ such that

\begin{align*}
(x^{\hat{\sigma}_1}y^{\hat{\sigma}_2})^k m_p m_{-u} &= x^\alpha y^\beta \in S \end{align*}

Let $\lfloor p_{\rho_i}\rfloor := \tilde{a}_{\rho_i}$ for all $i$ so that 

\begin{align*}
\lfloor p \rfloor &= (\lfloor p_{\rho_1} \rfloor, \dots, \lfloor p_{\rho_{\Sigma(1)}}\rfloor)\\
&= (\tilde{a}_{\rho_1}, \dots, \tilde{a}_{\rho_{|\Sigma(1)|}}) \end{align*}

and define the sets $A_+, A_-$ by 

\begin{align*}
m_p &= \frac{ \prod_{A_+}x_\rho^{a_\rho} \prod_{A_-} y_\rho^{a_\rho}}{\prod_{A_+} y_\rho^{a_\rho} \prod_{A_-} x_\rho^{a_\rho}} \end{align*} with $a_\rho = |\tilde{a_\rho}|$ for all $\rho$. Since $X_\Sigma$ is smooth, there exists $u_1\in L$ such that $u_1(e_\rho)=a_\rho$ for all $\rho \in (\sigma_1 \cap \sigma_2)(1) = \sigma_1(1) \cap \sigma_2(1)$. By the \textbf{Separation Lemma}, there exists $u_2\in L$ such that \[u_2\bigg|_{\sigma_1\cap \sigma_2}=0, u_2\bigg|_{\sigma_1\setminus \sigma_2}>0, u_2\bigg|_{\sigma_2\setminus \sigma_1}<0. \]

Since $x^{\hat{\sigma}_1}y^{\hat{\sigma}_2}$ is a generator of $I_{irr}$, up to torsion by $I_{irr}$ it suffices to consider only factors (over $S$) in the numerator and denominator of $m_p$ for variables $x_\rho$ and $y_\rho$ for $\rho \in \sigma_1\cup \sigma_2$. Now there exists $k_0 \in \mathbb{N}$ such that \\


\begin{align*}
(x^{\hat{\sigma}_1}y^{\hat{\sigma}_2})^{k_0} m_p &= \frac{ f_1}{\prod_{\sigma_2\cap A_+}y_\rho^{a_\rho} \prod_{\sigma_1 \cap A_-}x_\rho^{a_\rho}} \\
&= \frac{f_2}{\prod_{\sigma_1\cap\sigma_2\cap A_+}y_\rho^{a_\rho} \prod_{(\sigma_2\setminus \sigma_1) \cap A_+}y_\rho^{a_\rho}\prod_{\sigma_1\cap \sigma_2\cap A_-}x_\rho^{a_\rho}\prod_{(\sigma_1\setminus\sigma_2)\cap A_-}x_\rho^{a_\rho}} &[ \text{subdivision of sets }A_+, A_-] \end{align*}

for some monomials $f_1, f_2\in S$. Now $m_pm_{-u_1} = m_{p-u_1}$ since $u_1\in L$ gives that there exists $k_1\in \mathbb{N}$ such that 

\begin{align*}
(x^{\hat{\sigma}_1}y^{\hat{\sigma}_2})^{k_1} m_{p-u_1} &= \frac{ f_3}{\prod_{(\sigma_2\setminus \sigma_1)\cap A_+}y_\rho^{b_\rho}\prod_{(\sigma_1\setminus \sigma_2)\cap A_-}x_\rho^{b_\rho}} \end{align*}

for some $\{b_\rho\} \subset \N$ and some monomial $f_3\in S$, so that there exists $\ell\in \mathbb{N}$ such that $m_{p-u_1}m_{-\ell u_2} = m_{p-u_1-\ell u_2}$ since $\ell u_2\in L$ gives 

\begin{align*}
(x^{\hat{\sigma}_1}y^{\hat{\sigma}_2})^k m_{p-u_1-\ell u_2} &= f_4 \end{align*} for some monomial $f_4 = x^\alpha y^\beta\in S$, for some $\alpha,\beta\in \mathbb{N}^{|\Sigma(1)|}$.

\end{proof}

By Lemma 1, this implies there exists $k\in \mathbb{N}$ such that $(x^{\hat{\sigma}_1}y^{\hat{\sigma}_2})^k\cdot m_p \in M_{\Lambda(L)}$ for any choice of $\sigma_1, \sigma_2\in \Sigma(n)$ so that $coker(M_{\Lambda(L)} \hookrightarrow Im(f))=0$ for $\epsilon$ in the Effective Cone. \\[2cm]

In fact, this argument shows that for $p\in \mathcal{H}_L(0)$, $coker(M_{\Lambda(L)} \hookrightarrow Im(f))=0$ so that we have a resolution of the diagonal for any smooth toric variety $X_\Sigma$, without requiring the deformation $\mathcal{H}_L^\epsilon$. \\[2cm]

\section{Resolution of diagonal for toric D-M stack associated to global quotient of smooth toric variety by finite abelian group}

The resolution of the diagonal for $X_\Sigma$ smooth from the previous section also strengthens the construction of the diagonal object giving the kernel of a Fourier-Mukai transformation inducing the identity in $D^b_{Coh}(X_\Sigma')$ for the case that $X_\Sigma'=\faktor{X_\Sigma}{\mu}$ is a global quotient of a smooth toric variety $X_\Sigma$ by a finite abelian group $\mu$. \\

That is, since smooth toric varieties have an open cover by affine spaces: \[ X_\Sigma = \cup_{\sigma \preceq \Sigma(n)} U_\sigma \] with $U_\sigma = \text{Spec }\C[\sigma^\vee \cap M] \cong \A^m$ for $m=\text{rank }(M) = \text{rank } (N)$, we have that the same diagonal object \[ \pi(M_{\Lambda(L)}) = M_{\Lambda(L)} \otimes_{S[\Lambda(\tilde{L})]}S \] and corresponding resolution $(\mathcal{F}_{\mathcal{H}_L^\epsilon/\tilde{L}}, \partial)$ carry through in the case of a global quotient \[ X_\Sigma' = \faktor{X_\Sigma}{\mu}\] with $X_\Sigma$ smooth (not necessarily unimodular) and $\mu$ a finite abelian group. 

\newpage \Large \center \textbf{Appendix}
\appendix
\normalsize\flushleft
\textbf{Lemma 2 from Section~\ref{lemmata} for $X_\Sigma$ smooth, Fano, not necessarily unimodular}: For any $p \in \mathcal{H}_L^\epsilon(0)$ and $\sigma \in \Sigma(n)$, there exists $v_\sigma \in L$ with $(p+v_\sigma)_\rho <0$ only on $\sigma(1)$, and $(p+v_\sigma)_\rho >0$ for $\rho\not\in\sigma(1)$.  

\begin{proof} We previously noted that since $X_\Sigma$ is Fano, there exists $v\in L$ such that

\[ \begin{cases} v(u_\rho)<0 & \forall \text{ }\rho\in\sigma(1),\\ v(u_{\rho'}) \geq 0 & \forall \rho' \not\in\sigma(1) \end{cases} \] 

Now there exists $r\in \Z$ such that $(r\cdot v)_\rho <<0$ for all coordinates with $\rho\in\sigma(1)$ so that all $\rho\in \sigma(1)$ become much less than any coordinate $\rho'\not\in\sigma(1)$. Additionally, there exists a unique $\tilde{v}\in L$ such that \[ \begin{cases}  (p+\tilde{v})(u_\rho)= 0 & \forall\text{ }\rho\in\sigma(1)\\ (p+\tilde{v})(u_{\rho'}) > 0 & \forall\text{ }\rho'\not\in \sigma(1). \end{cases} \]

Now adding on sufficiently large multiples of $v$ and $\tilde{v}$ ensure that there exists $r, k\in \Z$ such that $v_\sigma = rv + k\tilde{v}$ gives that $p+v_\sigma$ satisfies

\[ \begin{cases} (p + v_\sigma)_\rho < 0 & \text{ iff } \rho\in\sigma(1),\\ (p+v_\sigma)_\rho  >0 & \text{ otherwise}. \end{cases} \] \end{proof}


\printbibliography

@article{coxhgscoordinatering,
  doi = {10.48550/ARXIV.ALG-GEOM/9210008},
  url = {https://arxiv.org/abs/alg-geom/9210008},
  author = {Cox, David A.},
  title = {Erratum to "The Homogeneous Coordinate Ring of a Toric Variety", along with the original paper},
  publisher = {arXiv},
  year = {1992},
  copyright = {arXiv.org perpetual, non-exclusive license}
}

@book{fulton,
 ISBN = {9780691000497},
 URL = {http://www.jstor.org/stable/j.ctt1b7x7vc},
 author = {William Fulton},
 publisher = {Princeton University Press},
 title = {Introduction to Toric Varieties. (AM-131)},
 urldate = {2022-08-09},
 year = {1993}
}

@book{C-L-S,
title = "Toric varieties",
author = "Cox, {David A.} and Little, {John B.} and Schenck, {Henry K.}",
year = "2011",
doi = "10.1090/gsm/124",
language = "English (US)",
isbn = "978-0-8218-4819-7",
volume = "124",
series = "Graduate Studies in Mathematics",
publisher = "American Mathematical Society",
address = "United States",
}

@misc{bayer-popescu-sturmfels,
  doi = {10.48550/ARXIV.MATH/9912247},
  
  url = {https://arxiv.org/abs/math/9912247},
  
  author = {Bayer, Dave and Popescu, Sorin and Sturmfels, Bernd},
  
  title = {Syzygies of Unimodular Lawrence Ideals},
  
  publisher = {arXiv},
  
  year = {1999},
  
  copyright = {Assumed arXiv.org perpetual, non-exclusive license to distribute this article for submissions made before January 2004}
}

@article{Beilinson1978,

author = {Beilinson, A. A.},
title = {Coherent Sheaves on Pn and Problems of Linear Algebra},
volume = {12},
journal = {Functional Analysis and Its Applications},
number = {3},
pages = {214 -- 216},
year = {1978},
doi = {10.1007/BF01681436},
url = {https://doi.org/10.1007/BF01681436}

}

@book{hartshorne,
  title={Algebraic Geometry},
  author={Hartshorne, R.},
  isbn={9780387902449},
  lccn={lc77001177},
  series={Graduate Texts in Mathematics},
  url={https://books.google.com/books?id=3rtX9t-nnvwC},
  year={1977},
  publisher={Springer}
}

@BOOK{miller-sturmfels,
    AUTHOR = "E.~Miller and B.~Sturmfels",
    TITLE = "Combinatorial Commutative Algebra",
    PUBLISHER = "Springer New York, NY",
    VOLUME = "1st Ed.",
    YEAR = 2005 }

@misc{bayer-sturmfels,
  doi = {10.48550/ARXIV.ALG-GEOM/9711023},
  url = {https://arxiv.org/abs/alg-geom/9711023},
  author = {Bayer, Dave and Sturmfels, Bernd},
  title = {Cellular Resolutions of Monomial Modules},
  publisher = {arXiv},
  year = {1997},
  copyright = {Assumed arXiv.org perpetual, non-exclusive license to distribute this article for submissions made before January 2004}
}

\end{document}